\documentclass[]{amsart}

\usepackage{amsmath}
\usepackage{amsthm}
\usepackage{amssymb}
\usepackage{graphicx}
\usepackage{amsfonts}
\usepackage{mathrsfs}
\usepackage[latin1]{inputenc}     
\usepackage[T1]{fontenc}          

\theoremstyle{plain}
\newtheorem{theorem}{Theorem}[section]

\theoremstyle{definition}

\newtheorem{remark}[theorem]{Remark}
\newtheorem*{remark*}{Remark}

\begin{document}

\title{Multidomain spectral method for Schr\"odinger equations}

\author[M.~Birem]{Mira Birem}

\author[C.~Klein]{Christian Klein}
\address
{Institut de Math\'ematiques de Bourgogne
9 avenue Alain Savary, BP 47870, 21078 Dijon Cedex, France}
\email{christian.klein@u-bourgogne.fr}

\begin{abstract}
A multidomain spectral method with 
compactified exterior domains combined 
with stable second and fourth order time integrators is 
presented for Schr\"odinger equations. 
The numerical approach allows high precision numerical studies of 
solutions on the whole real line. At examples for the linear and 
cubic nonlinear Schr\"odinger equation, this code is compared to 
transparent boundary conditions and perfectly matched layers 
approaches. The code can deal with 
asymptotically non vanishing solutions as the Peregrine breather being discussed 
as a model for rogue waves. It is shown that the Peregrine breather
can be numerically propagated with 
essentially machine precision, and that 
localized perturbations of this solution can be studied. 
\end{abstract}

\date{\today}

\subjclass[2000]{}

\thanks{We thank A.~Arnold, C.~Besse, A.~Boutet de Monvel, who interested us in 
the subject, J.~Frauendiener P.~Klein and K.~Roidot for helpful 
discussions and hints. }
\maketitle

\section{Introduction}
It is a common problem in the numerical treatment of partial 
differential equations (PDEs) that most of the PDEs appearing in 
applications are defined on $\mathbb{R}^{d}$ where $d$ is the 
dimension, i.e., on unbounded domains, whereas computers are 
generally able to deal with only finite domains. One popular approach 
in this context is to set up the problem on a finite domain and to 
impose \emph{transparent boundary conditions} (TBC) which  give the boundary data 
which would appear at 
this location if the initial value problem where treated on the whole 
$\mathbb{R}^{d}$, or\emph{ artificial boundary conditions} 
approximating TBC. To this end a certain relation between Dirichlet 
and von Neumann data at the boundary has to be satisfied called the 
Dirichlet to von Neumann map, see e.g.~\cite{hagstrom,semyon} 
for reviews of TBC 
approaches and in particular \cite{AntoineArnold} for Schr\"odinger equations. An alternative approach is 
given by \emph{perfectly matched layers} (PML) which were introduced by 
B\'erenger \cite{Berenger} for the Maxwell equations. The idea is to 
join absorbing layers at the boundaries of the computational 
domain. The to be solved PDE is deformed in the layers in a way that 
the solution in the computational domain is not affected, and that it 
is rapidly dissipated in the layers. This PML approach was adapted to 
various equations, see for instance \cite{nissan,note} for Schr\"odinger 
equations. We use here the approach \cite{Zhengpml} that is similar to \emph{exterior complex scaling} methods 
\cite{95}.

An interesting approach to treat PDEs on unbounded 
domains is to map infinite domains to finite domains. In a numerical 
context, this was apparently first used in \cite{orszag} and later for 
the Schr\"odinger equation in \cite{ladouceur}.
Whereas \cite{ladouceur} uses a finite difference method, 
\cite{orszag} incorporates one of the first applications of a 
spectral method with Chebyshev polynomials.   A very efficient 
variant of the latter mainly used in astrophysics, see for instance \cite{Lorene}, 
are multidomain methods with compactified exterior 
domains (CED). As in so called infinite element approaches, see for 
instance \cite{duque} for Schr\"odinger equations, the idea is to map 
the exterior domain with some M\"obius transformation to a finite 
domain. In this paper, we apply a multidomain spectral approach with 
CED to Schr\"odinger equations and compare this to TBC and PML 
approaches for concrete examples. The CED method clearly needs more 
resolution in space, but it will be shown that it is competitive 
especially in the nonlinear case since the TBC have to be implemented 
iteratively and since the parameters of the PML have to be optimized 
through several runs of the code for the studied example. In addition 
CED are attractive if high precision is needed since PML do not 
appear to provide sufficient accuracy in the NLS cases studied here, 
and since TBC are often not known and have to be implemented 
approximately. CED methods on the contrary are shown to provide the 
same accuracy in the linear and nonlinear cases.  Moreover the TBC 
and PML approaches here require compact support of the initial data 
in the computational domain whereas the CED approach just requires 
asymptotically bounded initial data and can even deal with an 
algebraic fall off of the solutions. 

In the focus of this paper are Schr\"odinger equations of the form
\begin{equation}
    i\partial_{t}u+\partial_{xx}u+V(|u|^{2},x)u=0
    \label{NLS},
\end{equation}
where $u=u(x,t):\mathbb{R}\times\mathbb{R}\mapsto \mathbb{C}$, and 
where the potential  $V$ is either a function of $x$ only or of $|u|^{2}$. In concrete 
examples, we study the case of the free linear Schr\"odinger equation 
with $V=0$ and of the focusing cubic nonlinear Schr\"odinger equation 
(NLS) with $V=2|u|^{2}$. It is straight forward to generalize the 
presented CED and PML approaches 
to more general potentials $V$, whereas this requires more work for 
TBC. The importance of the linear 
Schr\"odinger equation needs hardly be stressed since it is the 
governing PDE of quantum mechanics. In addition it arises in many 
applications as in quantum semiconductors \cite{28}, in electromagnetic 
wave propagation \cite{87}, in seismic migration \cite{33}, as the lowest order 
one-way approximation (paraxial wave equation) to the Helmholtz 
equation, under the notion of Fresnel equation in optics \cite{112}, or  in underwater acoustics \cite{129}.
But also the NLS equations appear 
in many applications in hydrodynamics, nonlinear optics and 
Bose-Einstein condensates mainly in the description of the modulation 
of waves. The cubic NLS which we study here in more detail was shown 
 to be completely integrable by Zakharov and Shabat \cite{ZS}. This 
implies that many explicit solutions to this equation are known. Of 
special interest in the context of rogue waves are \emph{rational 
breather} solutions to NLS, i.e., exact solutions with a slow 
algebraic fall off in both space and time. A popular candidate for 
rogue waves is  the Peregrine breather \cite{Peregrine}
\begin{equation}
    u_{Per} = \left(1-\frac{4(1+4it)}{1+4x^{2}+16t^{2}}\right)e^{2it}
    \label{peregrine},
\end{equation}
see  \cite{dubard} for generalizations thereof. It is claimed that this solution has been 
observed experimentally in rogue wave experiments in hydrodynamics \cite{cha1,cha2}, 
in plasma physics \cite{bailung} and in nonlinear optics 
\cite{kibler}. 

The Schr\"odinger equations are purely dispersive equations which 
means that the introduction of numerical dissipation should be 
limited if dispersive effects like rapid oscillations are to be 
studied. Spectral methods appear to be the ideal choice in this 
context because of their excellent approximation properties for 
smooth functions and because of the minimal introduction of numerical 
dissipation. If the solution is rapidly decreasing, Fourier spectral 
methods appear to be the most efficient tool, see for instance 
\cite{etna,klein} and references therein. If the initial data have compact 
support or are slowly decreasing as the Peregrine solution 
(\ref{peregrine}), Fourier methods are not ideal since the 
solution cannot be periodically continued as a smooth function even 
on large domains which leads to Gibbs phenomena at the boundaries. In this case multidomain spectral methods, for 
instance the Chebyshev collocation approach, see e.g.\
\cite{trefethen}, to be used in this paper appear to be better suited 
since they will show \emph{spectral convergence} if properly 
implemented, i.e., an exponential decrease of the numerical error 
with the number of collocation points for the approximation of analytic functions. 

If the Chebyshev collocation method is applied for several domains 
where the exterior ones are compactified (essentially $1/x$ is used 
as a coordinate there) and where appropriate matching conditions at 
the domain boundaries are implemented, a spectral method on the whole real line is 
obtained. It is straight forward to combine this approach with 
efficient and stable time integration schemes of high order. In 
this paper we present such an approach with the Crank-Nicolson scheme 
and an implicit fourth order Runge-Kutta method. It is shown that the 
CED method can be as easily implemented as PML. These approaches are 
tested at concrete examples for linear and nonlinear Schr\"odinger 
equations, where they are compared to TBC and PML methods known in the 
literature, for instance \cite{Zhengtbc,Zhengpml} for the cubic NLS. 
The Peregrine solution (\ref{peregrine}) is discussed as an example 
for slowly decreasing solutions. 

The paper is organized as follows:  we present in section 2 the 
multidomain spectral method with compactified exterior domain and how 
matching conditions are handled, in section 3 the used time 
integration schemes and how they are implemented, in section 4 the  
PML approach of \cite{Zhengpml}, and in section 5 the TBC approach for 
linear Schr\"odinger and cubic NLS. In section 6 we compare these 
approaches for Gaussian initial data for the linear Schr\"odinger 
equation, and in section 7 for the soliton of the cubic NLS equation. 
In section 8 we study with the CED approach the Peregrine breather and 
perturbations thereof. We add some concluding remarks in section 9. 

\section{Multidomain spectral method}\label{spectral}
In this section, we describe the  multidomain spectral method used to 
treat the spatial dependence of the 
function $u(x,t)$ in (\ref{NLS}). This means we divide the real line into a number of 
intervals, where we concentrate on the case of one finite interval and two infinite ones. Each of 
these intervals will be mapped to $[-1,1]$.  On the latter we use 
polynomial interpolation as in \cite{trefethen} to obtain a spectral 
approach. Since we are interested in solving 
(\ref{NLS}), a second order PDE in $x$, we impose at the boundaries 
of the intervals the matching condition that $u$ and its first 
derivative with respect to $x$ are continuous there. This ensures that 
we will get a spectral method on the whole real line.

\subsection{Lagrange polynomial and Chebyshev collocation method}

To approximate numerically the derivative of a function $f: 
[-1,1]\mapsto \mathbb{C}$, we use polynomial interpolation as 
presented for instance in \cite{trefethen}. For $l\in[-1,1]$, we 
introduce the $N+1$ \emph{Chebyshev collocation points} 
\begin{equation}
    l_{j}=\cos\left(\frac{j\pi}{N}\right),\quad j=0,\ldots,N
    \label{colpoints},
\end{equation}
where $N$ is some natural number.  We consider 
polynomial interpolation on these collocation points, i.e., 
the Lagrange polynomial $p(l)$ of order $N$ satisfying the relations
$p(l_{j})=f(l_{j})$, $j=0,\ldots,N$. The derivative of this 
polynomial is used as an approximation of the derivative of $f$ at 
the collocation points $l_{j}$, $f'(l_{j})\approx p'(l_{j})$. 
This is equivalent to approximate a 
derivative by the action of a \emph{differentiation matrix} $D$ on 
the vector $\mathrm{f}$ with the components $f_{j}:=f(l_{j})$, 
$j=0,\ldots,N$, i.e., $f'(l_{j})\approx (D \mathrm{f})_{j}$. The 
matrices $D$ are given in explicit form in 
\cite{trefethen}, a Matlab code to generate them can be found at 
\cite{trefethenweb}. Second derivatives of the function $f(l)$ will 
be approximated by $D^{2}\mathrm{f}$. 

This method is known to show spectral convergence, i.e., an 
exponential decrease of the numerical error  with the number $N$ of collocation points in approximating the 
derivative of an analytic function. In fact this 
\emph{pseudospectral} approach is equivalent to an approximation 
of a function $f$ by a truncated series of Chebyshev polynomials 
$T_{n}(l)$, $n=0,\ldots,N$, a \emph{Chebyshev collocation 
method}. This means we approximate $f$ as $f(l)\approx 
\sum_{n=0}^{N}a_{n}T_{n}(l)$, where the Chebyshev polynomials are 
defined as 
\begin{equation}
    T_{n}(l)=\cos(n\mbox{arccos}(l))
    \label{cheb},
\end{equation}
and where the 
\emph{spectral coefficients} $a_{n}$ are obtained by a collocation 
method, i.e., by imposing equality of $f$ and the sum of Chebyshev 
polynomials on the collocation points (\ref{colpoints}),
\begin{equation}
    f(l_{j})=\sum_{n=0}^{N}a_{n}T_{n}(l_{j}),\quad j=0,\ldots,N
    \label{coll}.
\end{equation}

Since the Chebyshev polynomials are related via (\ref{cheb}) to 
trigonometric functions, the coefficients $a_{n}$ in (\ref{coll}) 
can be computed via a 
\emph{fast cosine transformation} (fct) which is closely related to 
the \emph{fast Fourier transform} (fft), see the discussion in  
\cite{trefethen}. Since we use Matlab here and since the fct is in 
contrast to the fft not a precompiled command, it is much slower than 
the latter. Therefore we use here the pseudospectral approach via the 
Lagrange polynomial and corresponding differentiation matrices
in actual computations. But the relation to 
Fourier transforms implies that the well known fall off behavior of 
the Fourier coefficients of  an analytic function apply also to an 
expansion of a function on $[-1,1]$ in terms of  Chebyshev 
polynomials: for such a function $f$, the spectral 
coefficients $a_{n}$ in 
$f(l)=\sum_{n=0}^{\infty}a_{n}T_{n}(l)$ decrease exponentially 
with $n$ for $n\to\infty$. This implies that the numerical error in 
approximating an analytic function via a 
truncated series will also decrease exponentially. Since we will 
always consider functions in this paper where both real and imaginary 
part are  real analytic, we can test 
the consistence of our approach by studying the decrease of the spectral coefficients 
with $n$. This also allows to check whether sufficient spatial resolution 
is provided in each of the domains since the modulus of the 
coefficients should decrease to the wanted precision. Since we work 
here with double precision ($\approx 10^{-16}$), the maximal 
achievable precision 
is typically limited to  $\approx 10^{-14}$ because of rounding 
errors. We always aim at a resolution that the modulus of the 
coefficients decreases at least to $10^{-12}$ in each considered 
domain during the whole computation. 

An advantage of the  
use of a Chebyshev collocation method is that we can compute 
integrals as certain norms of the solution
conveniently with the \emph{Clenshaw-Curtis method}, see 
e.g.~\cite{trefethen}. The basic idea of this approach is as in 
(\ref{coll}) that the integrand is expanded in terms of Chebyshev 
polynomials, $$\int_{-1}^{1}f(l)dl\approx 
\sum_{n=0}^{N}a_{n}\int_{-1}^{1}T_{n}(l)dl=\sum_{n=0}^{N}w_{n}f(l_{n})$$ (the last 
step following from the collocation method (\ref{coll}) relating 
$a_{n}$ and $f(l_{n})$) where the 
$w_{n}$, $n=0,\ldots,N$ are some known weights (see 
\cite{trefethenweb} for a Matlab code to generate them). Since 
we sample $u$ already on Chebyshev collocation points, the 
integration is approximated by a scalar multiplication with the vector with 
components $w_{n}$. The integration scheme is also a spectral method, and the 
integral is thus computed with the same precision as the numerical 
solution to equation (\ref{NLS}).

An integration with respect to 
$x$ up to infinity (see the following subsection) implies for the integration with respect to 
$s=1/x$ division of 
the integrand by $s^{2}$ which vanishes for $x\to\infty$. Even if the 
integrand tends to zero sufficiently rapidly there, the numerical 
evaluation of expressions of the form `$0/0$' is problematic. 
In a spectral approach this division can be done in coefficient space 
which provides a numerically much more stable procedure. The approach 
is based on the well known recurrence formula for Chebyshev 
polynomials,
\begin{equation}
    T_{n+1}(l)+T_{n-1}(l)=2lT_{n}(x), \quad n=1,2,\ldots
    \label{chebrec}
\end{equation}
This formula allows to divide in coefficient space by $l\pm1$. We 
define for given Chebyshev coefficients $a_{n}$ coefficients $b_{n}$ 
via 
$\sum_{n=0}^{\infty}a_{n}T_{n}(l)=:\sum_{n=0}^{\infty}(l\pm1)b_{n}T_{n}(l)$. 
With (\ref{chebrec}) this implies the recursive relation 
\begin{align}
    a_{0} & =\pm b_{0}+\frac{1}{2}b_{1},
    \nonumber\\
    a_{1} & = b_{0}\pm b_{1}+\frac{1}{2}b_{2},
    \label{recur}\\
    a_{n} & = \frac{1}{2}b_{n-1}\pm b_{n}+\frac{1}{2}b_{n+1},\quad 
    n>1.
    \nonumber
\end{align}
Thus division with respect to $l\pm1$ can be done in coefficient 
space by solving (\ref{recur}) for the $b_{n}$. In practical 
applications, the series in $a_{n}$ and $b_{n}$ will be of course 
truncated. In addition numerical errors that 
$\sum_{n=0}^{N}a_{n}T_{n}(\pm1)$ does not vanish exactly. In this 
case the above procedure is applied to 
$\sum_{n=0}^{N}a_{n}(T_{n}(l)-T_{n}(\pm1))$.

\subsection{Multidomain method}
To use the spectral approach outlined in the previous subsection on 
the whole real line, 
we divide the latter into the three intervals 
$[-\infty,x_{l}]$, $[x_{l},x_{r}]$ and $[x_{r},\infty]$ 
($x_{l}<0$, $x_{r}> 0$) 
which will be denoted by I, II respectively III in the following. First 
we map the finite interval II to $[-1,1]$ with the linear 
transformation 
\begin{equation}
    x = x_{l}\frac{1+l}{2}+x_{r}\frac{1-l}{2},\quad l\in[-1,1]
    \label{xm}.
\end{equation}
Then we introduce Chebyshev collocation points 
$l_{j}^{II}$, $j=0,\ldots,N^{II}$, of the form (\ref{colpoints}), 
where $N^{II}\in\mathbb{N}$. The 
derivatives are approximated by the differentiation matrix $D^{II}$ 
obtained via polynomial approximation as explained in the previous 
subsection. 

On interval I, we use the mapping 
\begin{equation}
    x = \frac{2x_{l}}{1-l},\quad l\in[-1,1],
    \label{xl}
\end{equation}
which means we use a M\"obius transformation to map the infinite 
interval I to the compact interval $[-1,1]$. This approach is similar 
to infinite elements in finite elements methods, see \cite{duque} for 
Schr\"odinger equations. As before, $l$ is sampled on $N^{I}+1$, 
$N^{I}\in\mathbb{N}$ 
collocation points (\ref{colpoints}), and the corresponding 
differentiation matrix $D^{I}$ is introduced. 

Similarly we use on interval III the mapping 
\begin{equation}
    x = \frac{2x_{r}}{1+l},\quad l\in[-1,1],
    \label{xr}
\end{equation}
$N^{III}+1$, $N^{III}\in\mathbb{N}$,   collocation points (\ref{colpoints}) and the 
corresponding differentiation matrix $D^{III}$. Thus on each of the 
three intervals we have introduced a spectral method. Note that we 
could also use a single compactified domain $1/x_{l}<s<1/x_{r}$  where 
$s=1/x$. This produces gives the same results within numerical accuracy. We 
use here two compactified intervals to allow for different 
resolutions in these intervals since in the examples discussed, we 
consider waves propagating to the right and thus need more 
collocation points on 
this side of the axis. 

\begin{remark}
It is 
straight forward to introduce more than three intervals in the same 
way. This is of practical importance since it is known that the 
eigenvalues of the matrix $D^{2}$ grow as $N^{4}$, see the discussion 
in \cite{trefethen}. If $N$ is large, 
this can severely limit the achievable numerical precision in the 
solution of a PDE. Thus if 
high resolution is needed, it 
can be interesting to introduce an appropriate number of domains 
which are distributed in a way that the number $N$ of collocation 
points on each domain can be kept small though the spectral 
coefficients (\ref{coll}) decrease to machine precision. 
\end{remark}

\subsection{Boundary and matching conditions} \label{sectau}
At the boundaries between the different domains, the solution $u$ 
is required to be $C^{1}$. Since the Schr\"odinger 
equation is of second order in $x$, these conditions uniquely fix the 
solution for fixed $t$. Because the derivatives of $u$ are approximated 
with differentiation matrices, these conditions take with the 
spectral discretization introduced above the form for $x=x_{l}$
\begin{equation}
    u^{I}_{N^{I}}-u^{II}_{0}=0,\quad \frac{2}{x_{l}}
    \sum_{\alpha=0}^{N^{I}}D_{N^{I}\alpha}u^{I}_{\alpha}-
    \frac{2}{x_{l}-x_{r}}\sum_{\alpha=0}^{N^{II}}D_{0\alpha}u^{II}_{\alpha}=0
    \label{match1},
\end{equation}
and for $x=x_{r}$
\begin{equation}
    u^{II}_{N^{II}}-u^{III}_{0}=0,\quad 
    \frac{2}{x_{l}-x_{r}}\sum_{\alpha=0}^{N^{II}}D_{N^{II}\alpha}u^{II}_{\alpha}+
   \frac{2}{x_{r}}\sum_{\alpha=0}^{N^{III}}D_{0\alpha}u^{III}_{\alpha}=0
    \label{match2}.
\end{equation}
The first condition in both (\ref{match1}) and (\ref{match2}) ensures 
continuity of the function $u$ at the boundaries, the second the 
continuity of its derivative as approximated via differentiation 
matrices ($u'\approx Du$). 

These conditions will be implemented with the \emph{$\tau$-method} 
introduced by Lanczos \cite{tau}. The idea is to use the conditions 
(\ref{match1}) and (\ref{match2}) as additional equations for the 
equations obtained by discretizing the PDE (\ref{NLS}). We put 
$U=(u^{I},u^{II},u^{III})$, i.e., combine the vectors of the 
discretized solution in the three domains to a single vector. Thus 
equation (\ref{NLS}) is approximated by the system of ODEs in $t$ of 
the form
\begin{equation}
    i\partial_{t}U+\mathcal{L}U+VU=0,
    \label{NLSU}
\end{equation}
where $\mathcal{L}$ is the 
$(N^{I}+N^{II}+N^{III}+3)\times(N^{I}+N^{II}+N^{III}+3)$ matrix built 
from the blocks obtained from approximating the second derivative 
with respect to $x$ via differentiation matrices in the respective 
domains. As in \cite{trefethen}, we implement the matching conditions 
by replacing in the system (\ref{NLSU}) certain equations with 
(\ref{match1}) and (\ref{match2}): the line with the number 
$N^{I}$ (the count starts at zero) is replaced by the 
first condition in (\ref{match1}), line number $N^{I}+1$ by the 
second condition in (\ref{match1}), line number $N^{I}+N^{II}+1$ by 
the first condition in (\ref{match2}), line number $N^{I}+N^{II}+2$ 
by the second condition in (\ref{match2}); in all cases the 
corresponding component of right hand side is replaced by 0. 
This general approach will be slightly varied in accordance with the 
time integration schemes as detailed in the following section. 

It is known, see \cite{trefethen} and references therein, that the 
$\tau$-method does not  implement the boundary conditions exactly, 
but that the solution will satisfy the boundary conditions with the 
same spectral accuracy as the PDE. It turns out that it is 
numerically preferable to treat boundary conditions and PDE with the 
same approach. The alternative would be to solve the PDE with given 
boundary conditions (for instance vanishing conditions on the 
boundary) and use the homogeneous 
solutions (which depend on the used method for the time integration) 
to establish the matching conditions. One would have to deal with 
smaller matrices in this case. But since Matlab has very efficient 
algorithms for sparse matrices as the matrix 
$\mathcal{L}$, we use here this larger matrix. 
\begin{remark}\label{resolution}
    A spectral method means that a function is approximated on the 
    considered interval by globally smooth functions, here Chebyshev 
    polynomials. 
    The matching conditions that the function $u$ is $C^{1}$ at the 
    domain boundaries imply for the Schr\"odinger 
    equations being of second order in the spatial coordinate that a 
    smooth $u$ is obtained on the whole real line if the solution is 
    smooth on each domain. Consequently this approach leads to a 
    spectral method on the whole real line. Since 
    spectral methods are by nature global, this also implies that the 
    resolution as indicated by the Chebychev 
    coefficients\footnote{For a smooth function the modulus of the coefficients (\ref{coll}) 
    decreases for large $n$ and the order of magnitude of $|a_{n}|$ 
    for $n\sim N$ will be referred to as \emph{spatial resolution} in the 
    following.} in one domain affects the achievable accuracy in all 
    other domains. This means that if for instance the $|a_{n}|$ decrease in one 
    domain to $10^{-4}$ and in the others to $10^{-14}$, they will 
    decrease after a few time steps in all domains just to $10^{-4}$ 
    thus limiting the accuracy to this order of magnitude. 
    Consequently it is not possible to choose a small 
    resolution and thus lower accuracy in the exterior domains if one 
    is mainly interested in the finite domain. The lack of resolution 
    in the external domains will affect the achievable accuracy in 
    the finite domain. Note that the choice of the 
    same resolution in each domain does not imply that the number of 
    Chebychev coefficients is the same in all domains, but that the 
    the coefficients decrease to the same value in each domain which 
    can be archived by largely different values of $N^{I}$, $N^{II}$ 
    and $N^{III}$ depending on the problem and the choice of the 
    domain boundaries. 
\end{remark}

Note that we essentially use the coordinate $s=1/x$ in the 
compactified domain. With this coordinate, equation (\ref{NLS}) takes 
the form
\begin{equation}
        i\partial_{t}u+s^{4}\partial_{ss}u+2s^{3}\partial_{s}u+Vu=0
    \label{nlss},
\end{equation}
which is clearly singular at $s=0$ corresponding to $x\to\infty$. 
Because of this singular behavior, it is not necessary to impose a 
boundary condition at $s=0$. Regular solutions at infinity satisfy 
with (\ref{nlss}) $i\partial_{t}u(0,t)+Vu(0,t)=0$. However, it is 
possible to impose boundary values compatible with the 
 solution there  via a $\tau$-method, for instance the vanishing of the solution at 
 infinity for asymptotically vanishing solutions. 
This leads within numerical precision to the same solution on the 
whole line. But since 
we are also interested in solutions which do not vanish for 
$x\to\infty$ as the Peregrine solution (\ref{peregrine}), 
we do not impose conditions at infinity. Note that it 
is straight forward to implement (\ref{nlss}) with the 
differentiation matrices discussed above. All one has to do is to 
approach the operator $s^{4}\partial_{ss}+2s^{3}\partial_{s}$ which 
can be done with the same amount of work as on a finite domain.

\section{Time integration}
The spatial discretization in the previous section permits to 
approximate the PDE (\ref{NLS}) via a finite dimensional system of 
ordinary differential equations
(ODEs) (\ref{NLSU}). In this section we will present briefly the two schemes we 
employ to integrate this system, the Crank-Nicolson method and an 
implicit Runge-Kutta method of fourth order. Both are implicit 
because the high condition number (the ratio of the largest 
eigenvalue to the smallest) of the Chebyshev differentiation 
matrices makes explicit methods inefficient for stability reasons: 
since the largest eigenvalues of the matrix $D^{2}$ are of order 
$N^{4}$, the time steps would have to be chosen prohibitively small 
to satisfy stability criteria. 

\subsection{Crank-Nicolson}
The Crank-Nicolson (CN) scheme is an unconditionally stable method of 
second order which takes for the ODE $y'=f(y,t)$, $y\in 
\mathbb{R}^{g}$, $f:\mathbb{R}^{g}\mapsto\mathbb{R}^{g}$, the form
\begin{equation}
    y(t_{n+1})=y(t_{n})+\frac{h}{2}(f(y(t_{n}),t_{n})+f(y(t_{n+1}),t_{n+1}))
    \label{CN},
\end{equation}
where $h=t_{n+1}-t_{n}$ is the time step. Obviously the method is 
implicit for nonlinear $f(y)$.

For equation (\ref{NLSU}), the CN method yields
\begin{equation}
    (\hat{1}-ih\mathcal{L}/2)U(t_{n+1})=(\hat{1}+ih\mathcal{L}/2)U(t_{n})+
    \frac{ih}{2}\left(V(|U(t_{n+1})|^{2})U(t_{n+1})+V(|U(t_{n})|^{2})U(t_{n})\right)
    \label{CNsys}.
\end{equation}
The matching conditions (\ref{match1}) and (\ref{match2}) are implemented with a 
$\tau$-method as discussed in section \ref{sectau}: the same components of 
the matrix $\hat{1}-ih\mathcal{L}/2$ as in (\ref{NLSU}) are replaced by 
the matching conditions, the corresponding terms on the right hand 
side of the equation are put equal to zero. 
In the linear case ($V$ independent of $u$), the system (\ref{CNsys}) 
can be solved directly
numerically. In the nonlinear case, 
the  system is solved via a fixed point iteration which is stopped 
once the $L^{\infty}$ norm of consecutive iterates changes by less than $10^{-8}$.

As will be shown in the following sections, the CN method works well 
if not too high accuracy is needed. In practice a relative numerical error of
up to $10^{-5}$ is 
possible. If higher precision is wanted, a fourth order method is 
necessary. We mainly present CN here since the existing transparent 
boundary approaches for NLS use this method. Thus to compare the 
efficiency of TBC, PML and CED for the same time integration scheme, 
we implement CN in all cases. 

\subsection{Fourth order implicit Runge-Kutta method}
The general formulation of an $s$-stage Runge--Kutta method for the initial value problem
$y'=f(y,t),\,\,\,\,y(t_0)=y_0$ is as follows:
\begin{eqnarray}
 y_{n+1} = y_{n} + h      \underset{i=1}{\overset{s}{\sum}} \, 
 b_{i}K_{i}, \\
 K_{i} = f\left(t_{n}+c_ {i}h,\,y_{n}+h  
 \underset{j=1}{\overset{s}{\sum}} \, a_{ij}K_{j}\right),
 \label{K}
\end{eqnarray}
where $b_i,\,a_{ij},\,\,i,j=1,...,s$ are real numbers and
$c_i=   \underset{j=1}{\overset{s}{\sum}} \, a_{ij}$.

For the implicit Runge--Kutta scheme of order 4 
(IRK4) used here (Hammer-Hollingsworth method), one has
$c_{1}=\frac{1}{2}-\frac{\sqrt{3}}{6}$, 
$c_{2}=\frac{1}{2}+\frac{\sqrt{3}}{6}$, $a_{11}=a_{22}=1/4$,
$a_{12}=\frac{1}{4}-\frac{\sqrt{3}}{6}$, 
$a_{21}=\frac{1}{4}+\frac{\sqrt{3}}{6}$ and $b_{1}=b_{2}=1/2$. This 
scheme can also be seen as a 2-stage Gauss method. 

The system 
following from (\ref{K}) for (\ref{NLSU}) is written in the form 
\begin{align}
    (\hat{1}-iha_{11}\mathcal{L})K_{1} &
    = \Big((i\mathcal{L}U(t_{n})
    +iha_{12}\mathcal{L}K_{2}\nonumber\\
    &\left.+V\left(\left|U(t_{n})+h\sum_{j=1}^{2}a_{1j}K_{j}\right|^{2}\right)\left(U(t_{n})+h\sum_{j=1}^{2}a_{1j}K_{j}\right)\right)
    \nonumber,\\
    (\hat{1}-iha_{22}\mathcal{L})K_{2} & 
    =\Big(i\mathcal{L}U(t_{n})
    +iha_{21}\mathcal{L}K_{1}\nonumber\\
    &\left.+V\left(\left|U(t_{n})+h\sum_{j=1}^{2}a_{2j}K_{j}\right|^{2}\right)\left(U(t_{n})+h\sum_{j=1}^{2}a_{2j}K_{j}\right)\right)
    \label{Ksys}.
\end{align}
The matching conditions (\ref{match1}) and (\ref{match2}) are 
implemented via a $\tau$-method for the same indices as explained in 
section \ref{sectau}, but this time for both $K_{1}$ and $K_{2}$ 
for the matrices $\hat{1}-iha_{11}\mathcal{L}$ and 
$\hat{1}-iha_{22}\mathcal{L}$ respectively. The corresponding right-hand sides of the 
equations are put equal to zero. The resulting 
implicit system for $K_{1}$ and $K_{2}$ again requires an iterative 
solution in the nonlinear case. This is done by solving the equation 
in the form  (\ref{Ksys}) which gives some 
simplified Newton method. It shows rapid convergence in practice. 

In the linear case, the solution 
could be given in principle by inverting a 
$2(N^{I}+N^{II}+N^{III}+3)\times 2(N^{I}+N^{II}+N^{III}+3)$. However, 
this effective doubling of the dimension of the to be inverted matrix 
is not unproblematic because of the mentioned conditioning of the Chebyshev 
differentiation matrices. Therefore we iterate also in the linear 
case as explained above.

\section{Perfectly matched layers}
In this section we briefly summarize the PML approach as applied by 
Zheng \cite{Zhengpml} to Schr\"odinger equations. The basic idea is 
to replace the infinite intervals I and III of the CED approach  
by finite intervals of width $\delta>0$. In these intervals, the real 
coordinate $x$ is replaced by a complex coordinate which introduces 
into the purely dispersive Schr\"odinger equation some dissipation. 
The parameters of the deformation are to be chosen in a way that the 
solution $u$ is quickly damped in the layers without 
allowing reflections back to the computational domain II. In order 
for this concept to work, it is assumed that the initial data have 
compact support in zone II in contrast to the CED approach. 

The PML approach is thus  in the treatment of the 
spatial dependence of $u$ very similar to the CED approach, just that 
zone I is now given by $[x_{l}-\delta,x_{l}]$ and zone III by 
$[x_{r},x_{r}+\delta]$. The coordinate $x$ in (\ref{NLS}) is replaced 
by the coordinate $\tilde{x}$ defined as
\begin{equation}
    \tilde{x}=
    \begin{cases}
        x ,& x\in [x_{l},x_{r}] \\
        x+R\int_{x_{l}}^{x}\sigma(s)ds, & x\in[x_{l}-\delta,x_{l}]\\
        x+R\int_{x_{r}}^{x}\sigma(s)ds, & x\in[x_{r},x_{r}+\delta]	
    \end{cases}
    \label{xtilde},
\end{equation}
where $R=\exp(i\pi/4)$  and where $\sigma(s)$ is a  
positive damping function to be specified below. Assuming that the NLS 
equation is deformed in zones I and III to hold in the 
form (\ref{NLS}) with $x$ replaced by $\tilde{x}$, we get there
\begin{equation}
    i\partial_{t}u+\frac{1}{1+R\sigma(x)}\partial_{x}\left(
    \frac{1}{1+R\sigma(x)}\partial_{x}\right)u+Vu=0
    \label{NLSpml}.
\end{equation}

As in \cite{Zhengpml} we choose the damping function to be 
\begin{equation}
    \sigma(s)=
    \begin{cases}
        \sigma_{0}(s-x_{l})^{2}, & s\in[x_{l}-\delta,x_{l}] \\
         \sigma_{0}(s-x_{r})^{2},& s\in[x_{r},x_{r}+\delta]
    \end{cases}
    \label{sigma}
\end{equation}
where $\sigma_{0}$ is a to be specified positive parameter. Both 
$\delta$ and $\sigma_{0}$ determine the effective length 
$|\tilde{x}|$ on which the solution is dissipated. We always fix 
$\delta$ in the following experiments and vary $\sigma_{0}$ to 
identify optimal damping of the solution. For simplicity, the 
same parameters are chosen in both zones I and III.

The integration of the modified NLS equation (\ref{NLSpml}) is done 
as before. We introduce Chebyshev collocation points (\ref{coll}) in 
each of the three domains, in II as before for $x$, in I for 
$x$ in $\tilde{x}=R\sigma_{0}(x-x_{l})^{3}/3$ and similarly in zone 
III. Since the layer is chosen in a way that the solution is $C^{1}$ 
at the boundaries, we impose the same matching conditions 
(\ref{match1}) and (\ref{match2}). But since zones I and III are no 
longer unbounded, equation (\ref{NLSpml}) is not singular at the 
outer boundaries $x_{l}-\delta$ and $x_{r}+\delta$. At these points 
we impose the vanishing of $u$ as boundary condition which is again 
done with a $\tau$-method (see section \ref{sectau}). The time integration is performed as 
discussed in the previous section with the CN and IRK4 methods.

\section{Transparent boundary conditions}
In this section we present a brief summary of the TBC methods used 
here
for the free Schr\"odinger equation \cite{1} and the completely 
integrable NLS equation \cite{Zhengtbc}. For the spatial dependence in
both cases, we use the spectral method 
outlined in section 2 for the single domain $[x_{l},x_{r}]$. The time 
integration of the resulting system of ODEs is performed with the CN 
method (at least for NLS, we are not aware of other time integration 
schemes which have been applied so far).  The 
initial data are supposed to have compact support in $[x_{l},x_{r}]$, 
at least within numerical precision. For details of the 
approaches the reader is referred to 
\cite{1} and \cite{Zhengtbc}.

\subsection{Free Schr\"odinger equation}
We first consider the case of vanishing $V$ in (\ref{NLS}), i.e., the 
free Schr\"odinger equation,
\begin{equation}
    i\partial_{t}u+\partial_{xx}u=0,
    \label{schr}
\end{equation}
on a finite interval $[x_{l},x_{r}]$. For the 
spatial dependence we use the Chebyshev collocation method explained 
in section 2. Since we deal only with one domain in this section, we 
suppress the superscript $II$ here.

For the time integration we use the CN method, the boundary 
conditions are implemented with a $\tau$-method (see section \ref{sectau}). Recall that this means we replace the first line 
of the matrix $(\hat{1}-iD^{2}h/2)$ by $(1,0,\ldots,0)$ and the last 
line by $(0,\ldots,0,1)$ and the first and last component of the vector on the 
right hand side of (\ref{CN}) by $g_{0r}$ and $g_{0l}$ respectively, 
the Dirichlet boundary conditions. 
We call the resulting matrix on the left $\tilde{\mathcal{L}}$ and the vector 
on the right $(g_{0r},\tilde{u},g_{0l})^{t}$. Thus we get the solution
\begin{equation}
    u(t_{n+1})=\tilde{\mathcal{L}}^{-1}
    \begin{pmatrix}
        g_{0r} \\
        \tilde{u} \\
        g_{0l}
    \end{pmatrix}
    \label{CN2}.
\end{equation}

The transparent boundary conditions are found by effectively writing 
the Schr\"odinger equation in the form 
$(\exp(i\pi/4)\partial^{1/2}_{t}+\partial_{x})(\exp(i\pi/4)\partial^{1/2}_{t}-\partial_{x})u=0$ 
and by imposing an outgoing wave condition at the boundaries. This means that we 
have there the conditions
\begin{equation}
    \left.(\exp(i\pi/4)\partial^{1/2}_{t}u+\partial_{\mathbf{n}}u)\right|_{x=x_{l},x_{r}}=0
    \label{dvn},
\end{equation}
where $\partial_{\mathbf{n}}u$ denotes the normal derivative at the 
boundary pointing 
towards the exterior of the domain. Thus condition (\ref{dvn}) 
establishes the Dirichlet to von Neumann map for a solution to the 
Schr\"odinger equation, i.e., a relation between the Dirichlet data 
$g_{0l}$, $g_{0r}$ at the boundaries and the corresponding von 
Neumann data $g_{1l}$, $g_{1r}$ there. The fractional derivatives can 
be computed with different approaches, see \cite{6,1} for a CN 
discretization. Since in \cite{Zhengtbc} the latter approach by 
Antoine and 
Besse performed better, we concentrate on it which reads for the above setting
\begin{align}
     -g_{1l}(t_{n+1})& = F
     \sum_{k=0}^{n+1}\beta_{k}g_{0l}(t_{n+1-k})=:v_{l}+F\beta_{0}g_{0l}(t_{n+1}),
    \nonumber\\
     g_{1r}(t_{n+1})&= F
     \sum_{k=0}^{n+1}\beta_{k}g_{0r}(t_{n+1-k})=:v_{r}+F\beta_{0}g_{0r}(t_{n+1})
    \label{TBC},
\end{align}
 where $F=-\exp(-i\pi/4)\sqrt{2/h}$, and where
 \begin{equation}
     \beta_{k+2}=\beta_{k}\left(1-\frac{1}{k+1}\right), \quad k=0,1,\ldots,\quad  \beta_{0}=1, 
     \beta_{1}=-1
     \label{beta}.
 \end{equation}

On the other hand the normal derivatives in (\ref{dvn}) can be computed from the 
numerical solution (\ref{CN2}) by applying the differentiation matrix 
$D$. Thus we get the numerical approximation of the derivative, which 
we write at the boundary as 
\begin{align}
    g_{1r}(t_{n+1})&= \sum_{\alpha=0}^{N}D_{0\alpha}u_{\alpha}(t_{n+1})=: u_{r}(t_{n+1})+\gamma_{rr}g_{0r}(t_{n+1})+\gamma_{rl}g_{0l}(t_{n+1}),\nonumber\\
     g_{1l}(t_{n+1})&= \sum_{\alpha=0}^{N}D_{N\alpha}u_{\alpha}(t_{n+1})=:u_{l}(t_{n+1})+\gamma_{lr}g_{0r}(t_{n+1})+\gamma_{ll}g_{0l}(t_{n+1})   
    \label{ND}.
\end{align}
Equations (\ref{TBC}) and (\ref{ND}) imply (all quantities taken for 
$t=t_{n+1}$)
\begin{equation}
    \begin{pmatrix}
        \gamma_{rr}-F\beta_{0} & \gamma_{rl} \\
        \gamma_{lr} & \gamma_{ll}+F\beta_{0}
    \end{pmatrix}
    \begin{pmatrix}
        g_{0r} \\
        g_{0l}
    \end{pmatrix}=
    \begin{pmatrix}
        v_{r}-u_{r} \\
        -v_{l}-u_{l}
    \end{pmatrix}
    \label{rl}.
\end{equation}
In other words, the values for $g_{0l}$ and $g_{0r}$  resulting 
from (\ref{rl}) give transparent 
boundary conditions. 

Note that the above approach does not work in the presence of a 
non-constant potential $V(x)$ in the exterior of the interval 
$[x_{l},x_{r}]$. In this case $V(x)$ has to be either eliminated from 
the equation via a 
gauge transform, or it appears in the 
formal decomposition (\ref{dvn}) as part of a square root which has 
to be approximated, see for instance \cite{AntoineArnold,pauline} for 
the two approaches. 
The approach \cite{1} can be generalized in principle to higher order 
time integration schemes as the IRK4 method used here. To this end, 
a fourth order algebraic equation has to be solved, see also \cite{AntoineArnold}. We are not 
aware of any attempts in this direction. The advantage 
of CN is that the coefficients (\ref{beta}) can be given explicitly. 
The corresponding expressions for higher order schemes
are expected to be involved, and it is not clear whether 
they can be efficiently implemented.

\subsection{Cubic nonlinear Schr\"odinger equation}
It is well known that the cubic NLS equation, i.e., equation 
(\ref{NLS}) with $V=-2\rho |u|^{2}$,
\begin{equation}
    i\partial_{t}u+\partial_{xx}u-2\rho |u|^{2}u=0, \quad \rho=\pm1,\quad
    x\in\mathbb{R},\quad t>0
    \label{nls}
\end{equation}
is completely integrable, see \cite{ZS}; the equation is focusing for 
$\rho=-1$ and has solitonic solutions in this case, and is defocusing for $\rho=1$. 
The integrability of the equation implies 
that powerful solution techniques as Riemann-Hilbert problems exist. 
The latter could be used in \cite{anne} to establish the Dirichlet to 
von Neumann map for $x=const$. The found relations can be seen as 
above for the free linear Schr\"odinger equation as TBC for the 
equation. This was numerically implemented by Zheng in 
\cite{Zhengtbc}.  We 
give a brief summary of this approach. Denote the von Neumann data 
by $g_{1}(t)$ and the Dirichlet data by $g_{0}(t)$. 
In a first step one has to solve the advection type system for 
auxiliary quantities $L_{1,2}(t,s)$ and $M_{1,2}(t,s)$,
\begin{align}
    L_{1,t}-L_{1,s} & = ig_{1}(t)L_{2}+a(t)M_{1}+b(t)M_{2},
    \label{L1}  \\
    L_{2,t}+L_{2,s} & = -i\rho\bar{g}_{1}(t)L_{1}
    +\rho\bar{b}(t)M_{1}-a(t)M_{2},
    \label{L2}  \\
    M_{1,t}-M_{1,s} & = 2g_{0}(t)L_{2}+ig_{1}M_{2},
    \label{M1}  \\
    M_{2,t}+M_{2,s} & = 2\rho\bar{g}_{0}(t)L_{1}
	-i\rho\bar{g}_{1}(t)M_{1},
    \label{M2}
\end{align}
where 
\begin{equation}
    a(t) = i\rho \Im (g_{0}\bar{g}_{1}),\quad b(t) = 
    \frac{i}{2}(g_{0,t}+i\rho|g_{0}|^{2}g_{0}),
    \label{ab}
\end{equation}
with the  condition
\begin{equation}
    L_{1}(t,t) = \frac{i}{2}g_{1}(t),\quad M_{1}(t,t)=g_{0}(t),\quad 
    L_{2}(t,-t) = M_{2}(t,-t) = 0
    \label{cond}.
\end{equation}
With $M_{2}$ and $M_{1}$ given, one can then compute the 
Dirichlet to von Neumann map via
\begin{equation}
    g_{1}(t)=M_{2}(t,t)g_{0}(t)-e^{-\frac{i\pi}{4}}\left.\partial_{\tau}^{1/2}M_{1}(t,2\tau-t)
    \right|_{\tau=t}
    \label{dvnnls}.
\end{equation}
This provides a nonlinear variant of the 
Dirichlet-to von Neumann map (\ref{dvn}) in the linear case. 

To solve the system (\ref{M2}), it is convenient to introduce the characteristic coordinates 
\begin{equation}
    \xi = \frac{t+s}{2}, \quad \eta=\frac{t-s}{2}
    \label{char}.
\end{equation}
In the linearized case, the right 
hand sides in (\ref{L1})-(\ref{M2}) vanish approximately which implies
\begin{equation}
    L_{1}\sim \frac{i}{2}g_{1}(\eta),\quad M_{1}\sim g_{0}(\eta),\quad
    L_{2}\sim M_{2}\sim 0,
    \label{lin}
\end{equation}
which allows to recover the results from the previous subsection.

To solve this system numerically we introduce the following 
discretization: $t_{n}=nh$, $n=0,1,\ldots$ and 
$s^{n}=(-t_{n},-t_{n}+2h,\ldots,t_{n})$, i.e., the vector with the 
components $s^{n}_{m}=-t_{n}+2hm$, 
$m=0,1,\ldots,n$. The equations are integrated with 
the trapezoidal rule with respect to $\xi$ ($\eta$) whilst keeping 
$\eta$ ($\xi$) constant. Thus we obtain for the $\xi$ integration
\begin{equation}
    y^{n+1}_{m+1}-y^{n}_{m}=\frac{h}{2}(A^{n+1}x^{n+1}_{m+1}+
    A^{n}x^{n}_{m})
    \label{xiint},
\end{equation}
and for the $\eta$ integration
\begin{equation}
    y^{n+1}_{m}-y^{n}_{m}=\frac{h}{2}(A^{n+1}x^{n+1}_{m}+
    A^{n}x^{n}_{m})
    \label{etaint},
\end{equation}
where $y^{n}_{m}=y(t_{n},s_{m})$.

We first use (\ref{xiint}) and (\ref{M2}) to determine $M_{2}(t,t)$,
\begin{equation}
    M_{2,n+1}^{n+1}=M_{2,n}^{n}-i\rho h (a^{n+1}+a^{n})
    \label{M2end},
\end{equation}
and similarly for (\ref{L2})
\begin{equation}
    L_{2,n+1}^{n+1}=L_{2,n}^{n}+\frac{h}{2} \left(\frac{\rho}{2}(|g_{1}^{n+1}|^{2}
    +|g_{1}^{n}|^{2})+\rho(\bar{b}^{n+1}g_{0}^{n+1}+\bar{b}^{n}g_{0}^{n}) 
    -(a^{n+1}M_{2,n+1}^{n+1}+a^{n}M_{2,n}^{n})\right)
    \label{L2end}.
\end{equation}
With (\ref{etaint}) and (\ref{cond}) we get for (\ref{M1})
\begin{equation}
    M_{1,0}^{n+1}=M_{1,0}^{n}
    \label{M1start},
\end{equation}
and thus from (\ref{cond}) $M_{1,0}^{n}=0$ for all $n$, i.e., 
$M_{1}(t,-t)=0$. This implies  $L_{1}(t,-t)=0$. Thus we have the 
functions $L_{i}$, $M_{i}$, $i=1,2$ for $m=0$ and $m=n$.

For $t_{n+1}$ and $1\leq m\leq n$, we get for the system (\ref{L1}) 
to (\ref{M2}) with (\ref{xiint}) and (\ref{etaint})
\begin{align}
    L_{1,m}^{n+1} & =L_{1,m}^{n}+\frac{h}{2}\left( 
    ig_{1}^{n+1}L_{2,m}^{n+1}+a^{n+1}M_{1,m}^{n+1}
    +b^{n+1}M_{2,m}^{n+1}+ig_{1}^{n}L_{2,m}^{n}+a^{n}M_{1,m}^{n}
    +b^{n}M_{2,m}^{n}\right),
    \label{L1d}  \\
    L_{2,m}^{n+1} & =L_{2,m-1}^{n}+\frac{h}{2}
    \left(-i\rho\bar{g}_{1}^{n+1}L_{1,m}^{n+1}
    +\rho\bar{b}^{n+1}M_{1,m}^{n+1}-a^{n+1}M_{2,m}^{n+1}
    \right.\nonumber\\
    &\left.
    -i\rho\bar{g}_{1}^{n}L_{1,m-1}^{n}
	+\rho\bar{b}^{n}M_{1,m-1}^{n}
	-a^{n}M_{2,m-1}^{n}\right),
    \label{L2d}  \\
    M_{1,m}^{n+1} & =M_{1,m}^{n}+\frac{h}{2}\left( 2g_{0}^{n+1}
    L_{2,m}^{n+1}+ig_{1}^{n+1}M_{2,m}^{n+1} +2g_{0}^{n}
	L_{2,m}^{n}+ig_{1}^{n}M_{2,m}^{n}\right),
    \label{M1d}  \\
    M_{2,m}^{n+1} & =M_{2,m-1}^{n}+\frac{h}{2}\left( 
    2\rho\bar{g}_{0}^{n+1}L_{1,m}^{n+1}
	-i\rho\bar{g}_{1}^{n+1}M_{1,m}^{n+1}
	+2\rho\bar{g}_{0}^{n}L_{1,m-1}^{n}
	-i\rho\bar{g}_{1}^{n}M_{1,m-1}^{n}\right).
    \label{M2d}
\end{align}
Defining for fixed $n,m$ the vector $X^{\gamma}$ a with components 
$L_{1},L_{2},M_{1},M_{2}$ where the index $\gamma$ refers for 
$1,2,3,4$ to $L_{1},L_{2},M_{1},M_{2}$ respectively, we can write this 
system in the form
\begin{equation}
    \sum_{\delta=1}^{4}\mathcal{A}_{\gamma\delta}X^{\delta}=V^{\gamma},\quad 
    \gamma=1,2,3,4
    \label{formal},
\end{equation}
where $A$ and $V$ follow from (\ref{L1d}) to (\ref{M2d}). Obviously 
the left-hand sides 
of  (\ref{L1d}) to (\ref{M2d}) have the same form for all $m$, 
whereas the right-hand sides follow by 
shifting indices. Thus the system can be directly solved by inverting 
the $4\times4$ matrix $\mathcal{A}$ for given $n$ and $m$. 

With $M_{2}$ and $M_{1}$ given, one can then compute the 
Dirichlet to von Neumann map (\ref{dvnnls}) with the Antoine-Besse approach as in 
the previous subsection. 
Here we essentially use the approach by Zheng \cite{Zhengtbc}  with the 
following changes: for the spatial discretization we use the 
pseudospectral method of section 2. This system is solved with 
CN by a fixed point iteration. In each step of 
the iteration, we compute $u(x_{l},t)$ and $u(x_{r},t)$ as well as 
the normal derivatives via the differentiation matrix $D$. The 
found normal derivatives via the procedure outlined in this subsection is 
then used as von Neumann data for the next step in the iteration. 
Thus in contrast to Zheng, we iterate the Dirichlet to von Neumann 
maps for $x_{l}$ and $x_{r}$ and the CN relations at the same time.

\section{Numerical study of the free Schr\"odinger equation}
\begin{figure}[htb!]
   \includegraphics[width=0.49\textwidth]{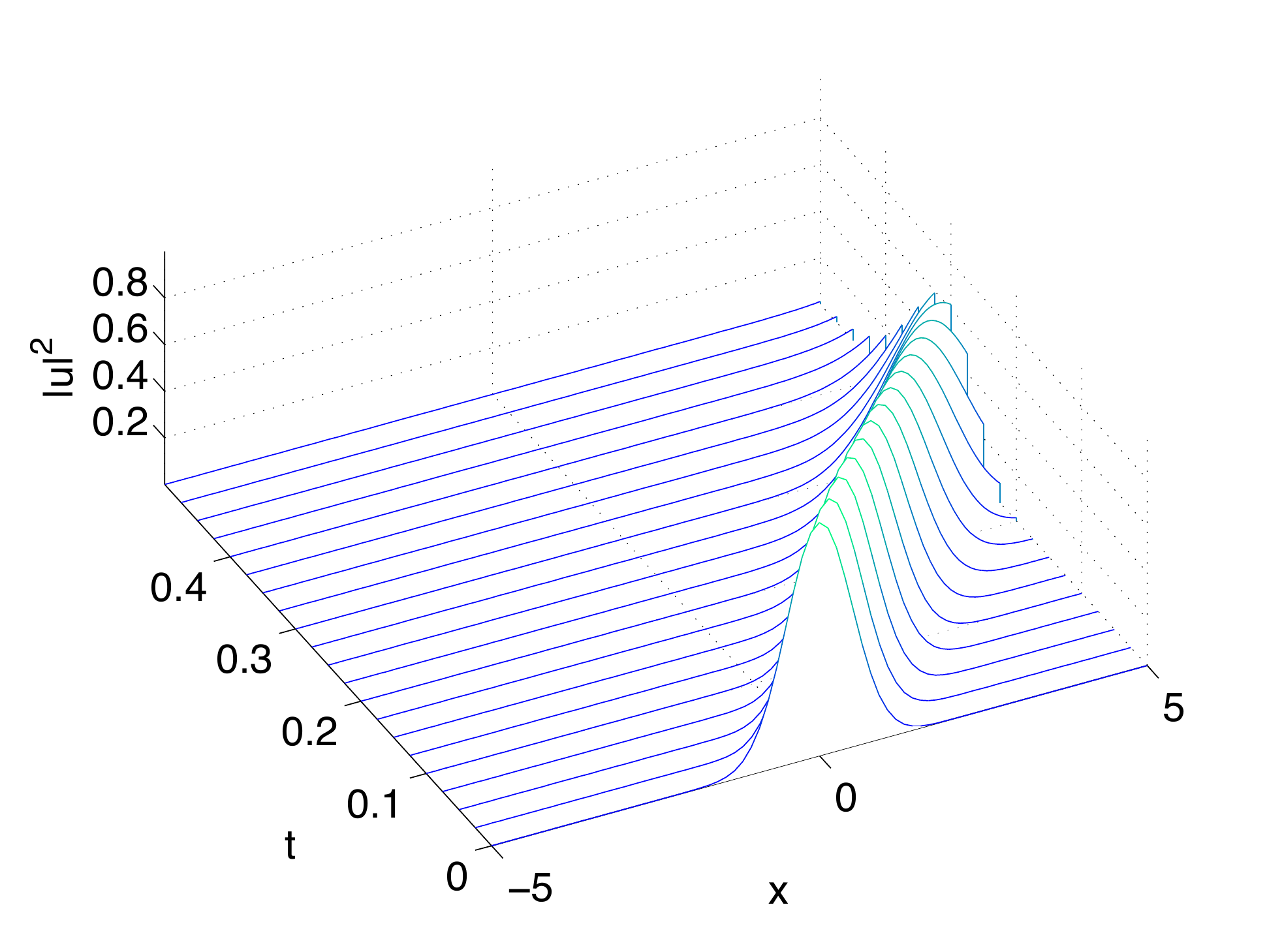}
   \includegraphics[width=0.49\textwidth]{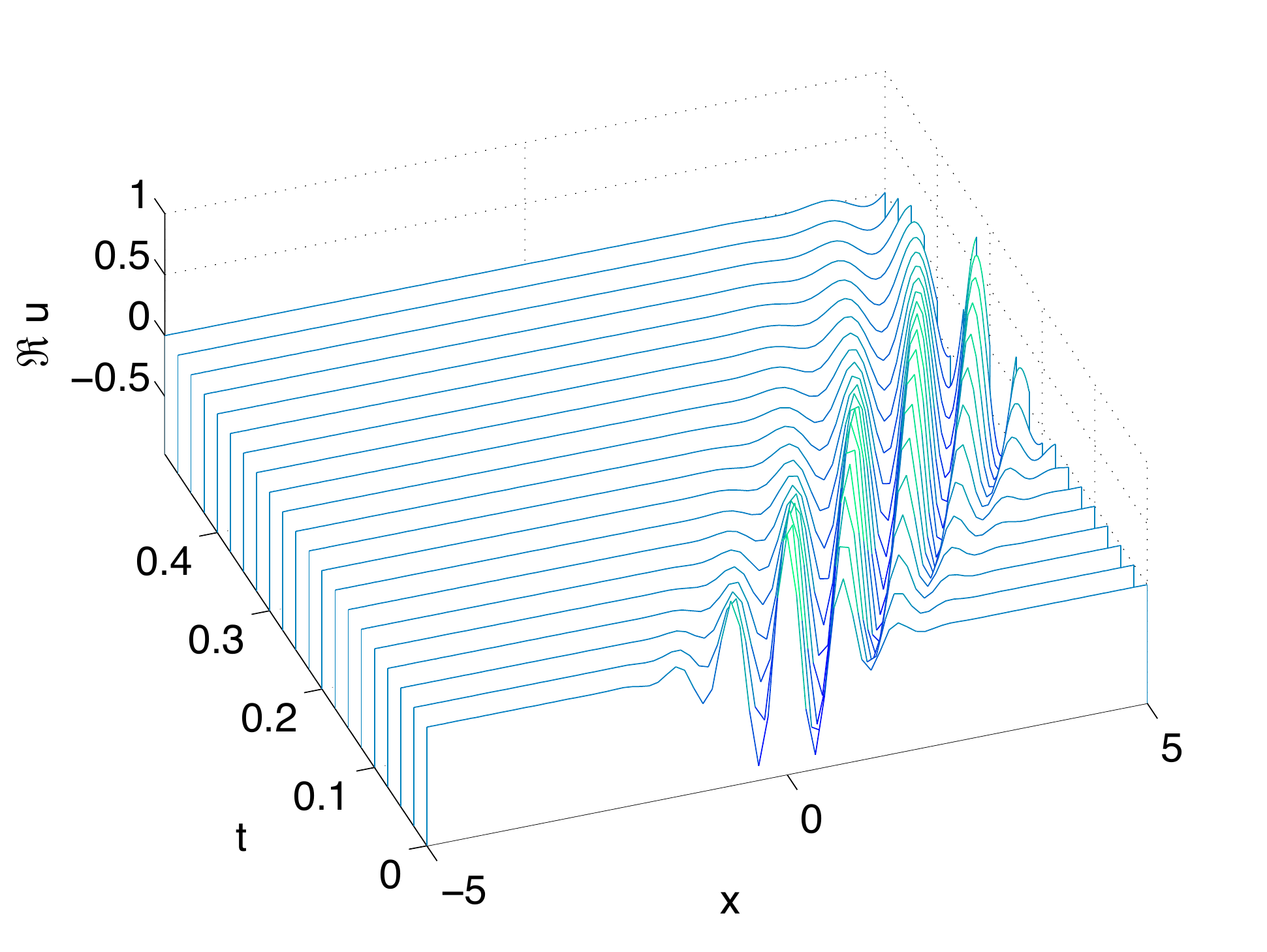}
 \caption{Exact solution (\ref{schrex}) to the free Schr\"odinger equation 
  being used as a test of the various numerical approaches; on the left the modulus squared of the solution, on the 
 right its real part. }
 \label{schroex}
\end{figure}

In this section, we will compare the CED method with PML and TBC of 
the previous sections for the free Schr\"odinger equation.  As 
in \cite{Zhengpml}, we consider as an example the exact solution
\begin{equation}
    u(x,t) =     
    \frac{1}{\sqrt{1+4it}}\exp\left(-\frac{x^{2}-8ix+8^{2}it}{1+4it}\right)
    \label{schrex}
\end{equation}
shown in Fig.~\ref{schroex} and give $u(x,0)$ as initial data for the to be tested numerical 
approaches. The numerical solution is then compared for $t\leq 1/2$ with 
the exact solution. It is shown that all methods are able to produce 
in principle the same accuracy in the linear case, but that they are 
not equally efficient: whereas TBC needs the least spatial 
resolution, higher order time 
integration schemes have not been explored so far. The latter can be easily done for CED, but for this 
method  to be of high precision, the Chebyshev coefficients have to 
decrease to the same order of magnitude
on the whole real line. The PML approach is a good compromise 
in this sense, since high order time integration can be easily used, 
and since only comparatively low spatial resolution is needed in the 
layers. But for this advantage is in practice more than offset by the need to 
optimize the parameters $\sigma_{0}$, $\delta$ through several runs 
of the code for the same initial data. 

It can be seen in Fig.~\ref{schroex} that the modulus of the solution 
(\ref{schroex}) is being slowly dispersed whilst the maximum travels 
with the constant speed $c=16$ to the right and hits the right boundary of 
domain II at $t=0.3125$. The real part of the solution in the same 
figure shows the typical oscillations of solutions of dispersive 
equations which also have to be accurately resolved 
numerically. 

The solution (\ref{schrex}) can be constructed from the well known 
solution to the Schr\"odinger equation for Gaussian initial data by 
exploiting the  {\it Galilei 
invariance} of the Schr\"odinger equations of the form 
(\ref{NLS}). This means, that if $u(t,x)$ is a solution, 
then so is 
\begin{equation}\label{eq:galilei}
\hat{u}(x,t)=u(x-ct,t)e^{ic x/2 -ic^{2}t/4},
\end{equation}
with $c\in \mathbb{R}$ some finite speed.
\begin{remark}
    It is convenient in practical computations for data with a single 
    maximum to exploit the Galilean invariance of the equation by 
    going to a reference frame essentially travelling with the maximum. 
    In this way, it is straight forward to choose an optimal 
    resolution for different domains in CED approaches. This can be 
    even done by going to an accelerated frame and thus changing the 
    form of the 
    to be solved equation as for instance in \cite{KP2013}. For the 
    example (\ref{schrex}), this would imply the study of the 
    dispersed Gaussian with a maximum fixed at the origin. Since the 
    solution is rapidly decreasing, Fourier techniques would 
    certainly be most efficient here.  However, 
    the goal of the present paper is to compare the efficiency of 
    various numerical approaches in 
    dealing with a maximum passing through the domain boundary which 
    is why we do not switch to a moving frame here and use the 
    example of \cite{Zhengpml}. 
\end{remark}

\begin{figure}[htb!]
   \includegraphics[width=0.49\textwidth]{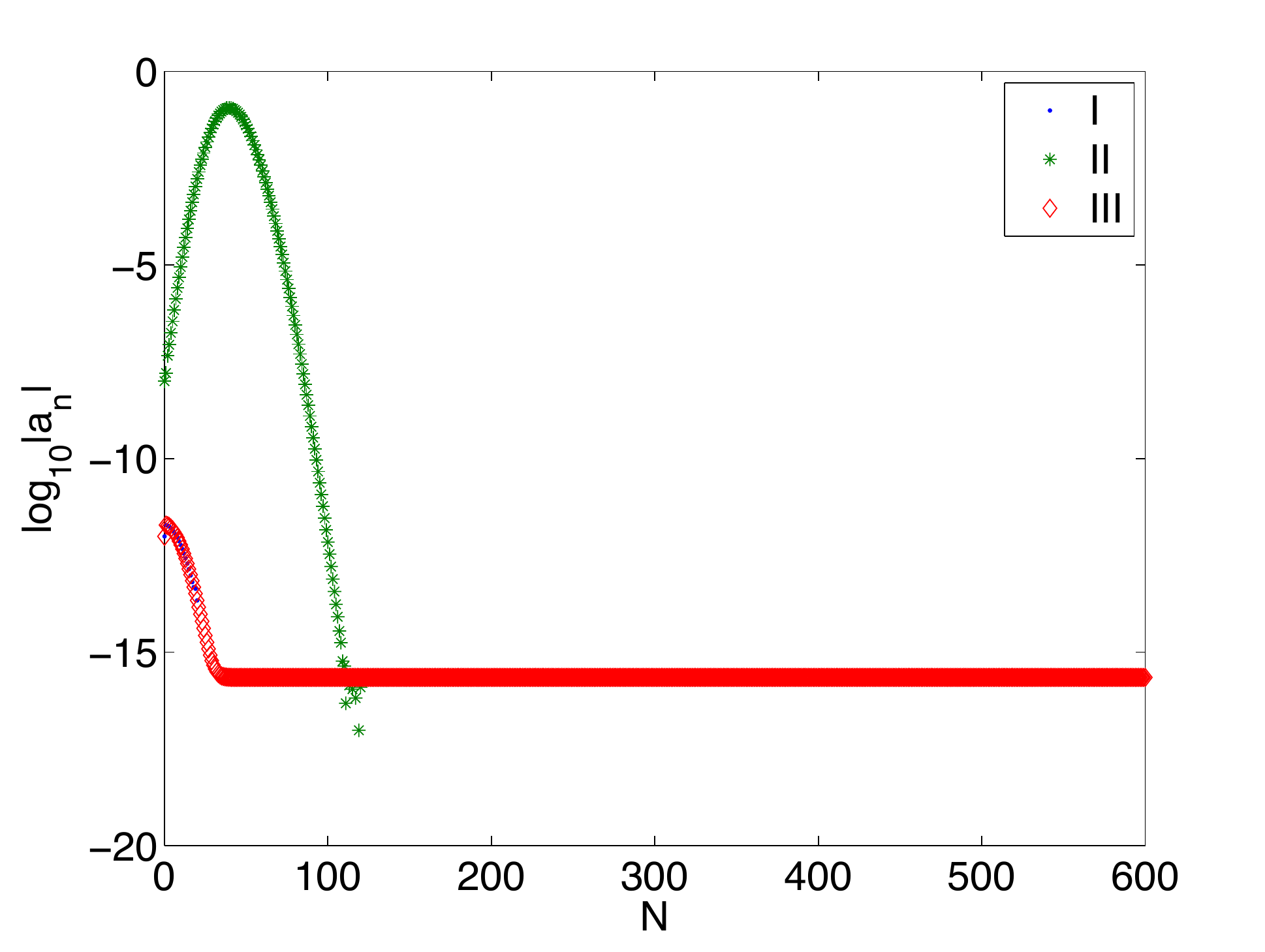}
   \includegraphics[width=0.49\textwidth]{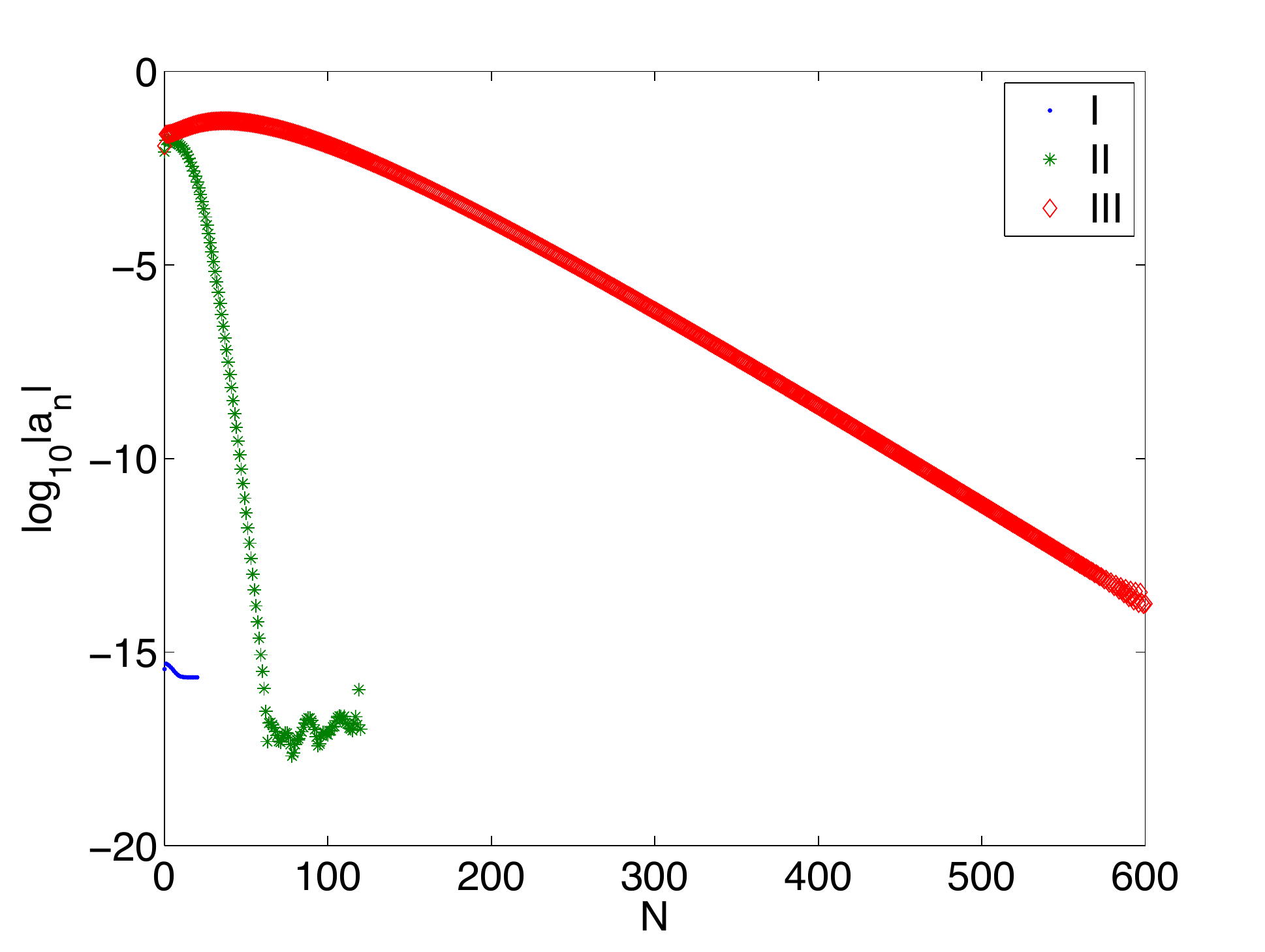}
 \caption{Chebyshev coefficients (\ref{coll}) for the  solution 
 (\ref{schrex}) in the three domains; on the left at  $t=0$, on the 
 right for $t=0.5$.}
 \label{schroeres}
\end{figure}

For the CED method, we choose $N^{I}=20$, $N^{II}=120$ and 
$N^{III}=600$ Chebyshev points in the respective domains. It can be 
seen in Fig.~\ref{schroeres} that the Chebyshev coefficients  
(\ref{coll}) decrease in each of the domains to machine precision, 
both for the initial and the final time of the computation. Note that 
the Chebyshev coefficients in domain I are always completely of the order of 
machine precision (they are covered in the left figure in 
Fig.~\ref{schroeres} by the coefficients in domain III) which means 
one could just impose vanishing Dirichlet condition at $x_{l}$. The 
domain is mainly kept here to illustrate the approach. 

As in \cite{Zhengpml}, we will discuss in the following the numerical 
error $\Delta$ as the $L^{2}$ norm of the difference between 
numerical solution $u$ and exact solution $u_{ex}$ normalized by the 
$L^{2}$ norm of the exact solution,
\begin{equation}
    \Delta=\frac{||u-u_{ex}||_{2}}{||u_{ex}||_{2}}.
    \label{Delta}
\end{equation}
For PML and TBC, both are 
computed  for $x\in[x_{l},x_{r}]$, for CED for $x\in\mathbb{R}$. This 
is done with the Clenshaw-Curtis algorithm outlined in section 
\ref{spectral}.

\begin{figure}[htb!]
   \includegraphics[width=0.49\textwidth]{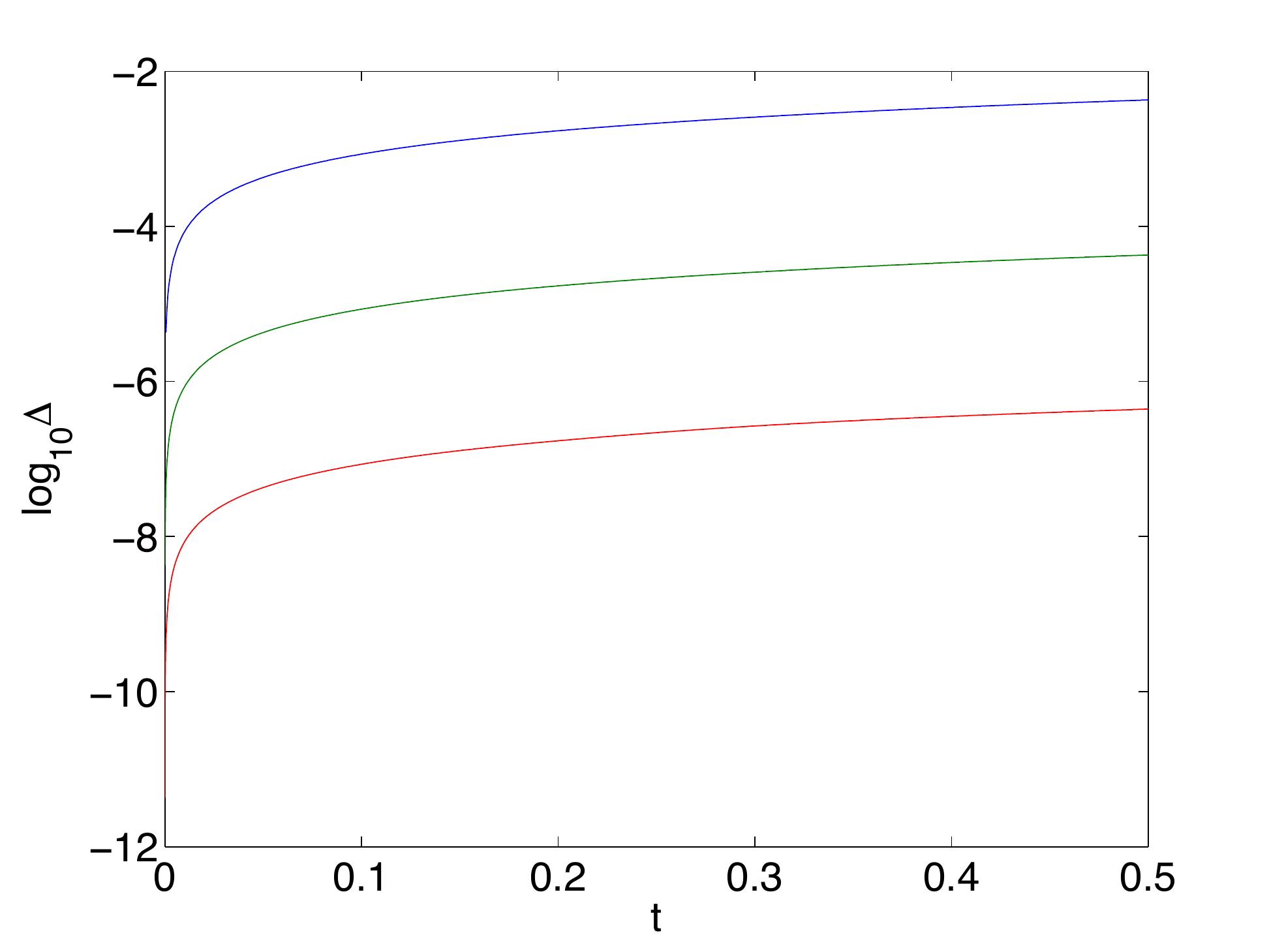}
   \includegraphics[width=0.49\textwidth]{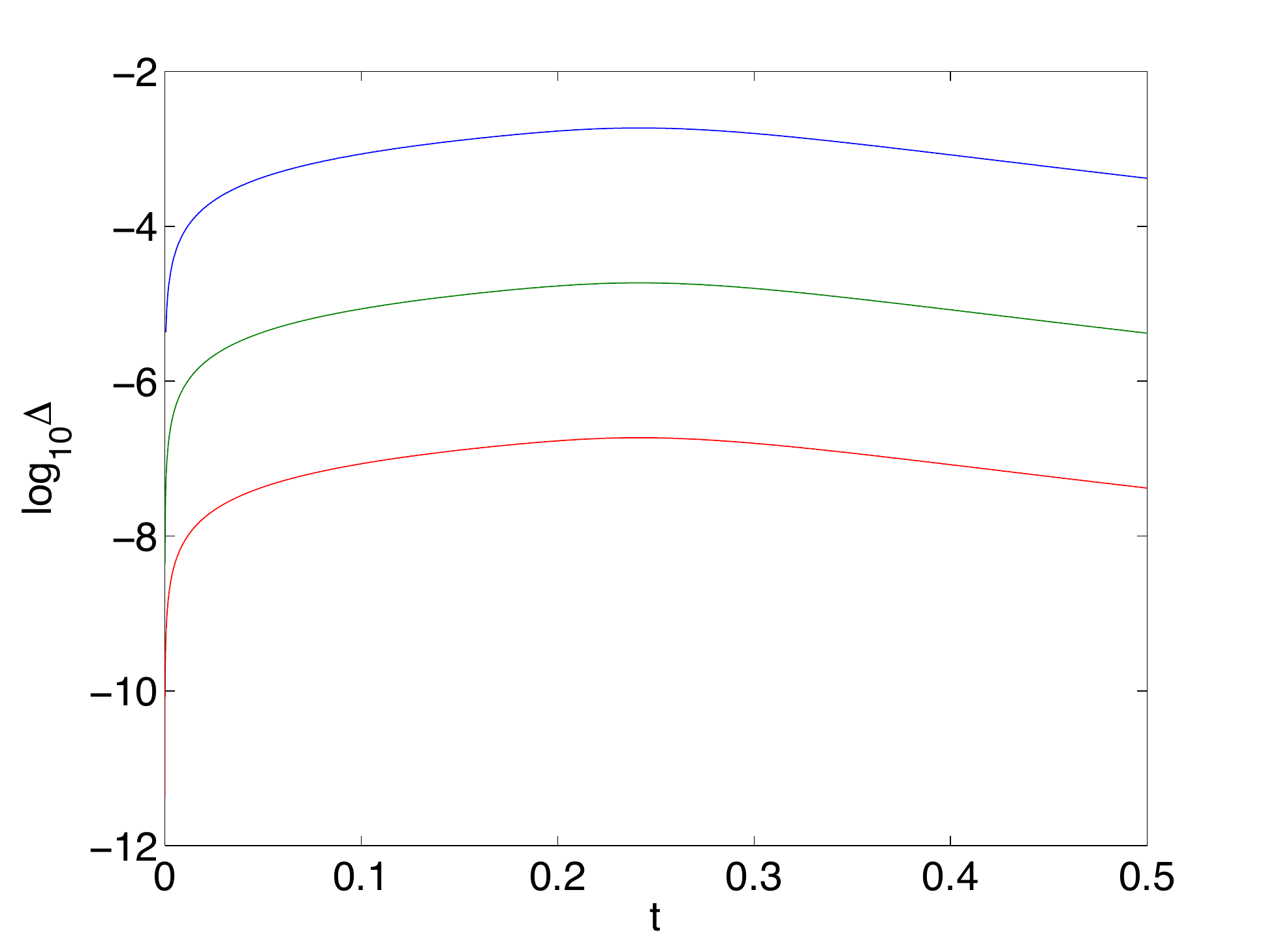}
 \caption{Normalized $L^{2}$ norm $\Delta$ of the difference between numerical 
 and exact solution (\ref{Delta}) with CN method for $N_{t}=10^{3}$, 
 $N_{t}=10^{4}$ and $N_{t}=10^{5}$ time steps (from top to bottom); 
 on the left for the CED method the error on the whole line, on the 
 right for the TBC method in the computational domain.}
 \label{schroemd}
\end{figure}

We first compare the various approaches with second order CN method in 
all cases. The numerical error $\Delta$ for the CED method can be seen 
in Fig.~\ref{schroemd} on the left. Because of the matching 
conditions at the domain boundaries, this is a global method on the 
whole real line, and no effect should be visible in this approach 
when the maximum crosses the boundary at $t=0.3125$. It can be seen 
that this is indeed the case, the numerical error increases in time 
due to a piling up of the errors in the time integration. In 
Fig.~\ref{schroemd}, the numerical error is shown for three 
resolutions ($N_{t}=10^{3},10^{4},10^{5}$). It can be clearly seen 
that the CN method shows the expected second order behavior. A 
precision of the order of $10^{-7}$ can be reached in this example, 
for higher accuracies a fourth order method as discussed below would 
be clearly the better choice. 

In Fig.~\ref{schroemd}, the same situation is studied on the right 
for the TBC approach. Here we only have domain II where we use the 
same resolution as before ($N^{II}=120$). Since both CED and TBC can 
be treated explicitly, the latter is considerably more efficient 
since only one domain is considered. The numerical error behaves 
slightly differently in this case as can be recognized in 
Fig.~\ref{schroemd}: once considerable parts of the mass of the 
solution leave around $t\approx 0.3$ the domain, the numerical error 
there decreases slightly since the solution itself is almost zero in 
part of the domain. There are no spurious 
reflections at the boundary if sufficient resolution in time is 
provided, the conditions are really transparent 
within numerical precision which is theoretically expected, see 
\cite{AntoineArnold}. The CN scheme again shows the expected 
second order behavior.

\begin{figure}[htb!]
   \includegraphics[width=0.49\textwidth]{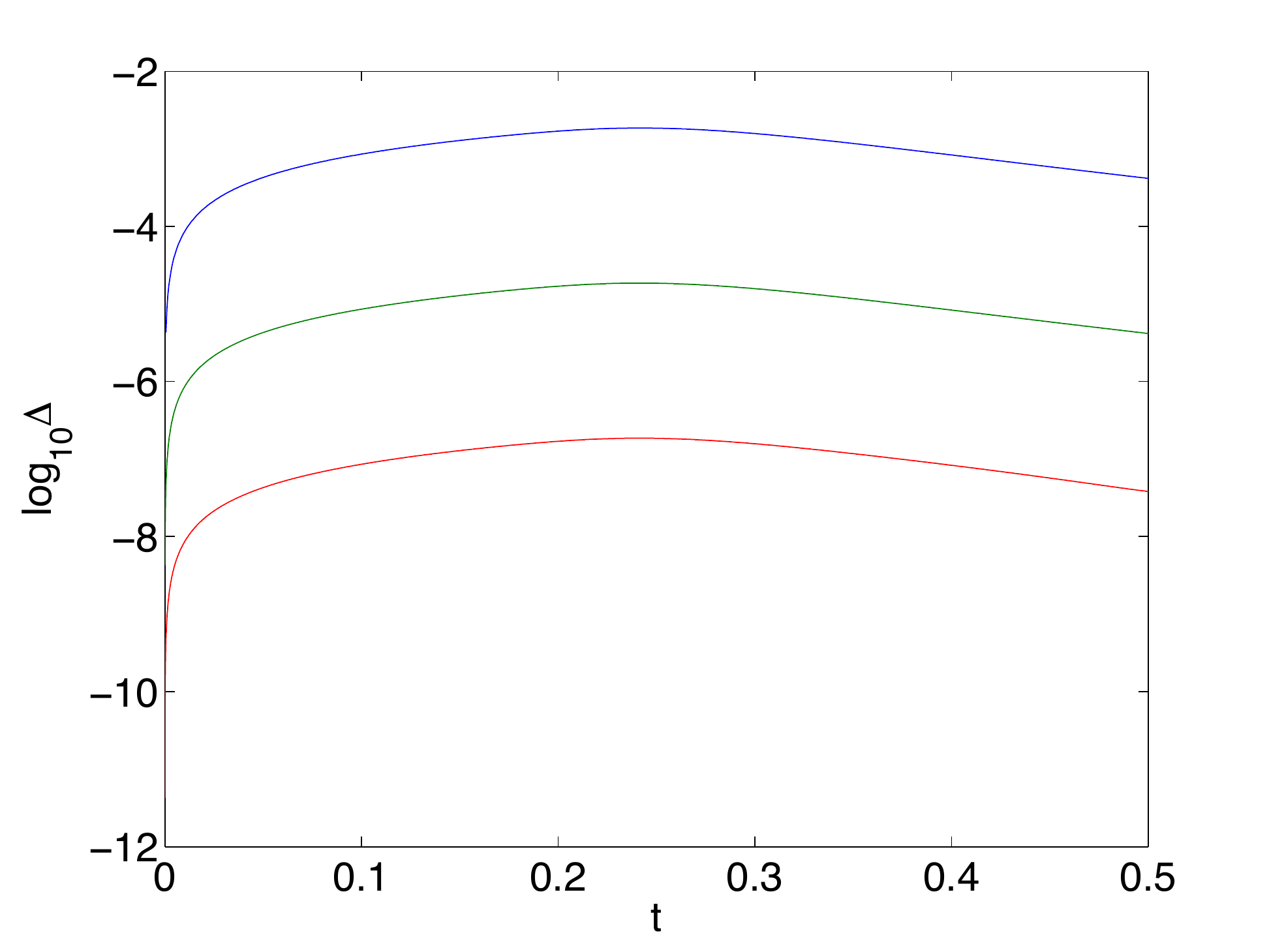}
   \includegraphics[width=0.49\textwidth]{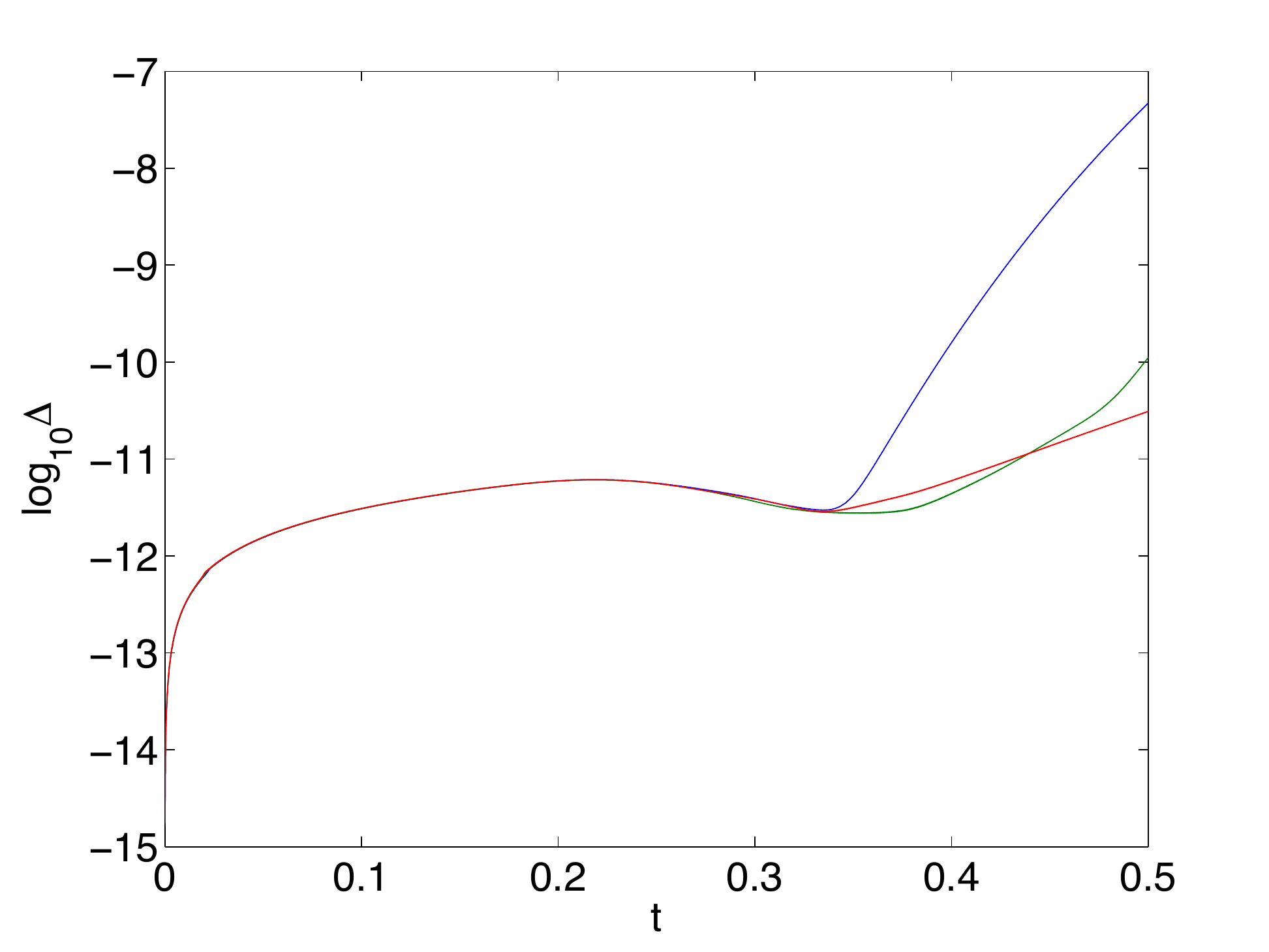}
 \caption{On the left, normalized $L^{2}$ norm $\Delta$ of the difference between numerical 
 and exact solution (\ref{Delta}) with CN method for $N_{t}=10^{3}$, 
 $N_{t}=10^{4}$ and $N_{t}=10^{5}$ time steps (from top to bottom) 
 for the PML method with $\sigma_{0}=50$ and $\delta=0.5$; on the right the PML approach with IRK4 method 
 and $N_{t}=10^{4}$ for various values of $\sigma_{0}$ in 
 (\ref{sigma}) ($\sigma_{0}=40,50,60$ corresponding to the blue, 
 green, and red curve respectively).}
 \label{schroepml}
\end{figure}

As can be seen in Fig.~\ref{schroepml} on the left, the numerical 
error in the PML case with CN is very similar to  the 
TBC case. Again second order convergence of the time integration 
scheme is observed, and the numerical error (\ref{Delta}) decreases 
once most of the mass leaves the computational domain. We use here 
$N^{I}=20$ and $N^{III}=50$ Chebyshev collocation points and find 
that the Chebyshev coefficients (\ref{coll}) decrease to $10^{-12}$ in the layers
throughout the computation. This 
shows that much less spatial resolution is needed here in domain III 
compared to the CED approach. Thus for the shown accuracy range, the 
PML method is less efficient than TBC, but more so than CED.

\begin{figure}[htb!]
   \includegraphics[width=0.49\textwidth]{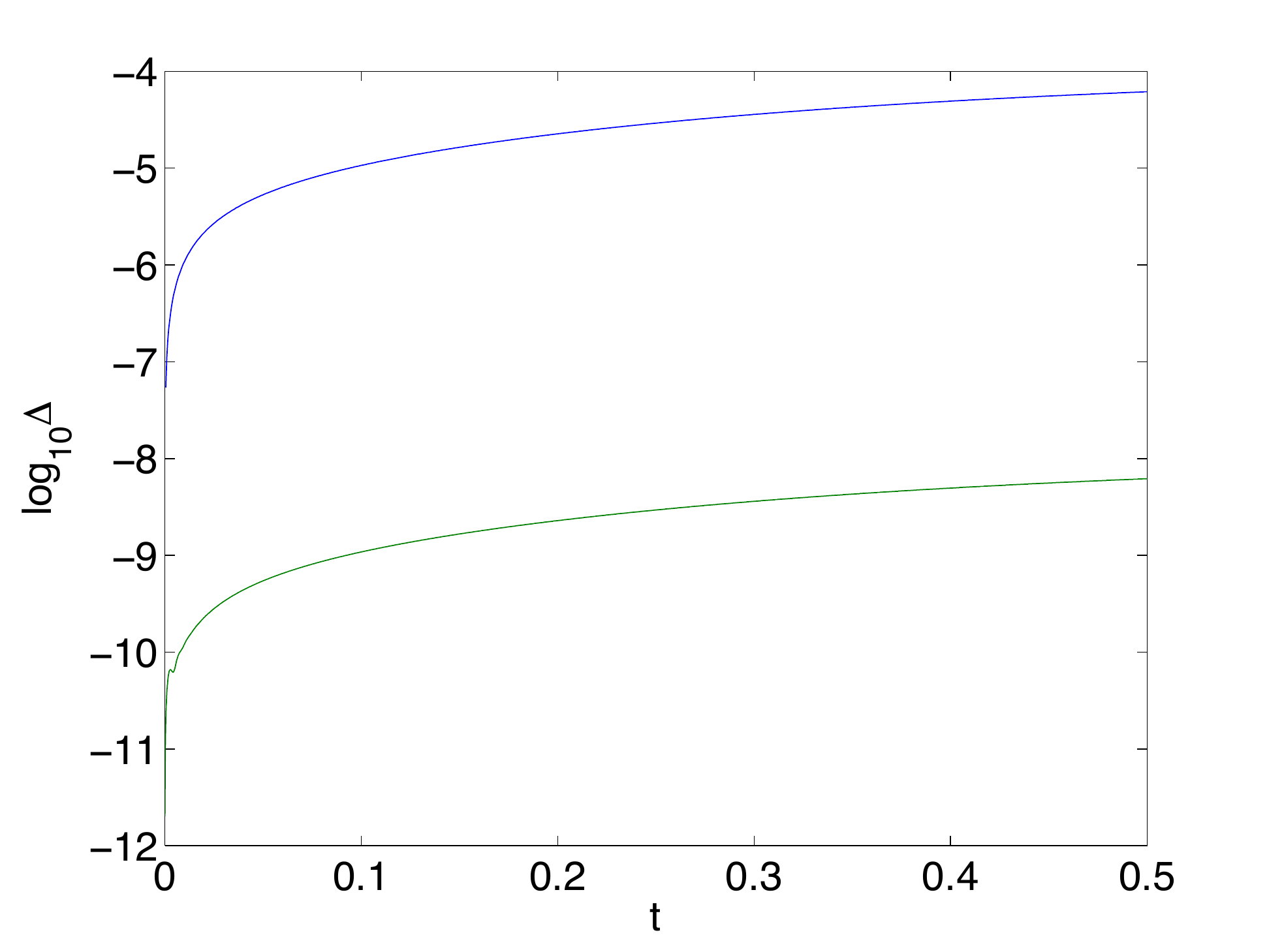}
   \includegraphics[width=0.49\textwidth]{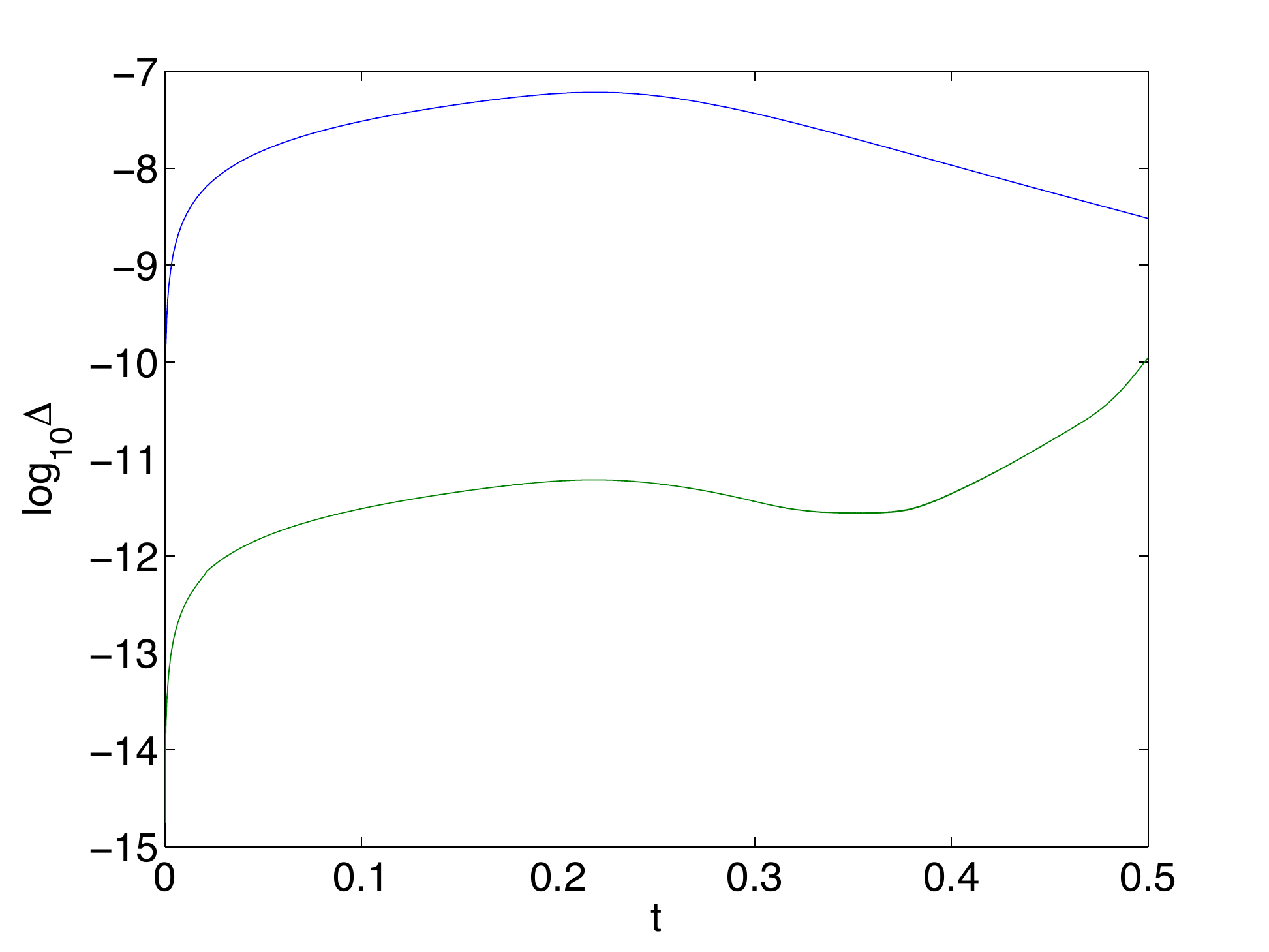}
 \caption{Normalized $L^{2}$ norm $\Delta$ of the difference between numerical 
 and exact solution (\ref{Delta}) with IRK4 method for $N_{t}=10^{3}$ 
 and  
 $N_{t}=10^{4}$ time steps (from top to bottom); 
 on the left for the CED method the error on the whole line, on the 
 right for the PML method with $\delta=0.5$ and $\sigma_{0}=50$ in the computational domain.}
 \label{schroeirk4}
\end{figure}

The above examples with the second order scheme CN in 
Fig.~\ref{schroemd} and \ref{schroepml} show that it is 
hardly possible to reach machine precision with such an approach. 
Therefore we apply the fourth order IRK4 method here to illustrate 
the difference. It is an advantage of both PML and CED that it is 
straight forward to combine them with convenient time integration 
schemes. It can be seen in Fig.~\ref{schroeirk4} that essentially the 
same accuracy can be reached for $N_{t}=10^{3}$ as with $10^{5}$ time 
steps for CN in Fig.~\ref{schroemd}, and that with $N_{t}=10^{4}$ one 
essentially reaches machine precision for the CED method. 

For the PML method, the same situation is shown on the right of 
Fig.~\ref{schroeirk4}. The behavior for $N_{t}=10^{3}$ is as expected 
from the CED case and the CN approach for PML in Fig.~\ref{schroepml}. 
But for the even higher resolution $N_{t}=10^{4}$, an effect of the 
right layer can be seen. After the maximum of the solution hits the 
domain boundary at $t=0.3125$, the numerical error increases again in 
the computational domain due to reflections from the layer. This is 
why the values for $\delta$ and $\sigma_{0}$ in (\ref{sigma}) were 
fixed after studying the numerical error with the IRK4 method in 
Fig.~\ref{schroepml} on the right. Since both $\delta$ and 
$\sigma_{0}$ have the same effect to increase the effective length 
$|\tilde{x}|$ of the layer where the dissipation is active, we fix as 
in \cite{Zhengpml} $\delta=0.5$ and vary $\sigma_{0}$ to minimize the 
numerical error. It can be seen in Fig.~\ref{schroepml} that the 
numerical error $\Delta$ increases strongly for $\sigma_{0}=40$ once 
the maximum of the solution enters the layer. This is much less the 
case for larger $\sigma_{0}$, but $\sigma_{0}=60$ produces a larger 
error near $t=0.3125$ than $\sigma_{0}=50$. Since it is most 
important to resolve the solution well in the computational domain 
whilst there is still considerable mass there, we chose 
$\sigma_{0}=50$ for the computations. The need to optimize the value 
of $\sigma_{0}$ essentially by trial and error, i.e., by additional 
runs of the code for the same initial data, is a clear disadvantage 
of the PML approach. Note however that in the present example, there is no 
visible effect by choosing either of these three values for 
$\sigma_{0}$ with the CN method on the left of Fig.~\ref{schroepml} 
since the accuracy of the time integration scheme is too low. 

Note 
also that 
the error shown for CED is obtained on the whole axis whereas it is for PML 
only for the computational domain. Thus it is not surprising that the 
error is smaller for PML at later times since the solution is close 
to zero after a certain time in the computational domain, whereas the 
maximum of the solutions is being tracked outside this domain in the 
CED approach. 

\section{Numerical study of the NLS soliton}
\begin{figure}[htb!]
   \includegraphics[width=0.49\textwidth]{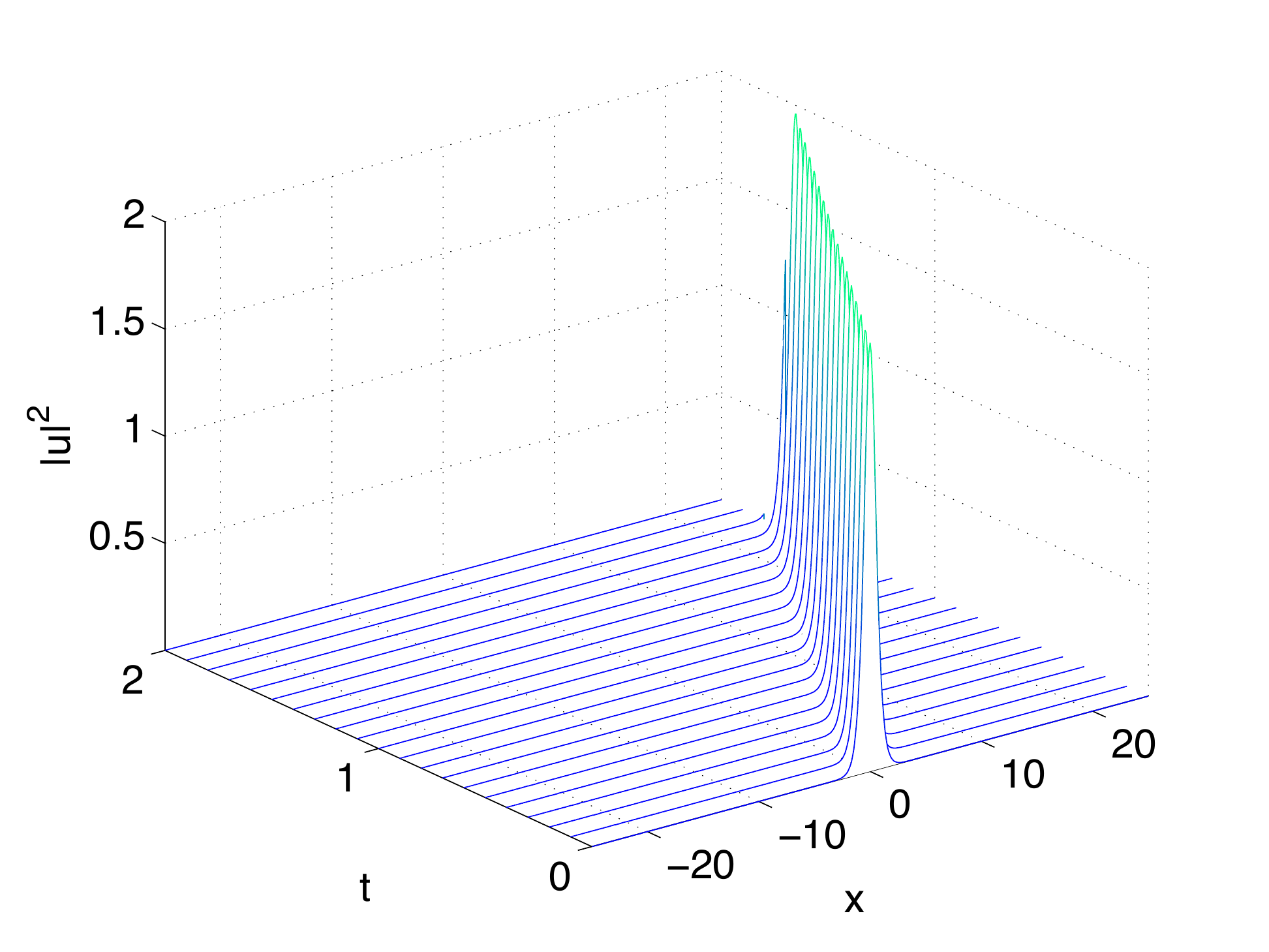}
   \includegraphics[width=0.49\textwidth]{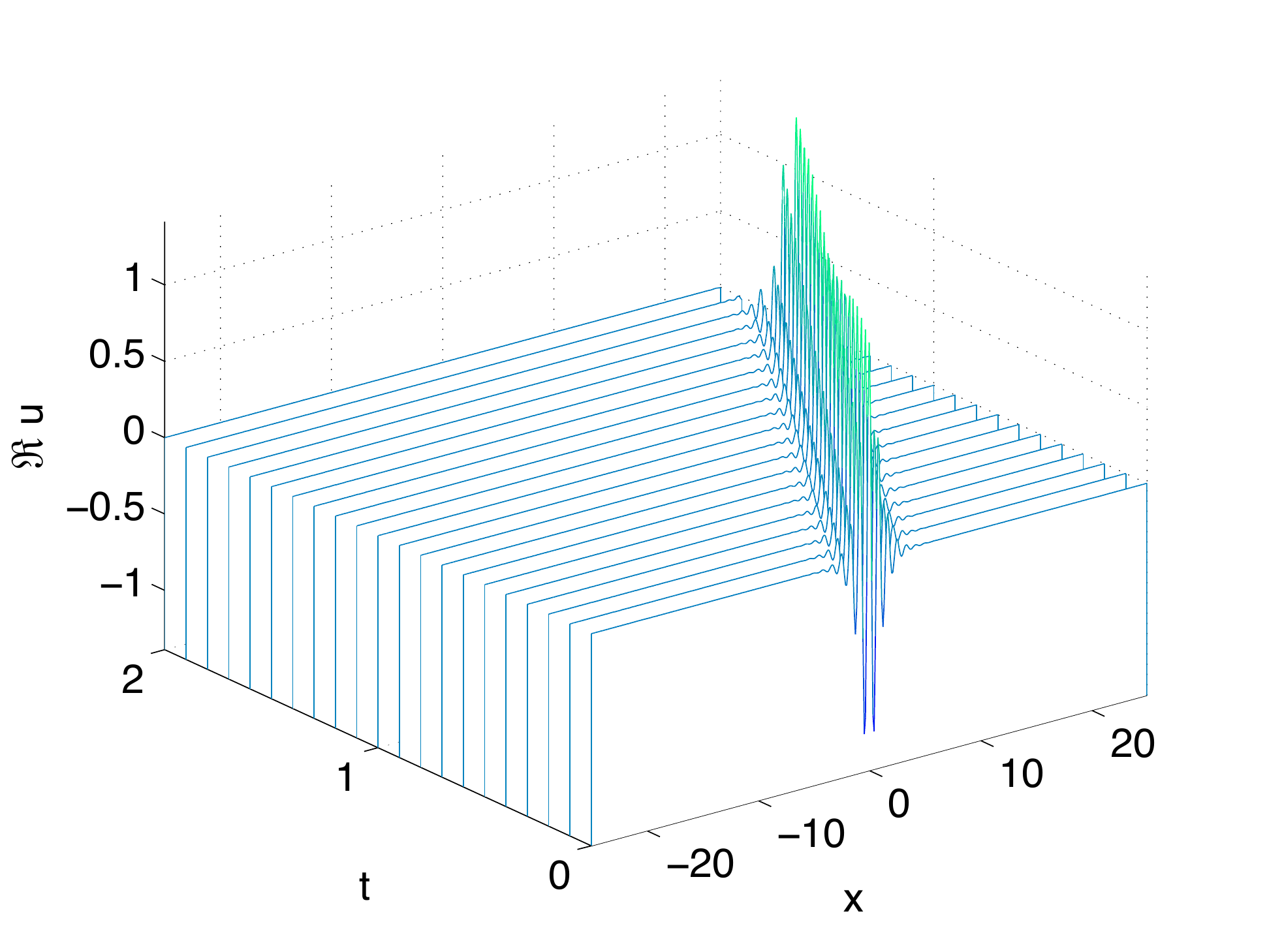}
 \caption{Soliton solution (\ref{sol}) to the cubic NLS equation 
being used as a test for the various numerical approaches; on the left the modulus squared of the solution, on the 
 right its real part. }
 \label{nlssol}
\end{figure}

In this section we study the CED, TBC and PML approaches for the NLS 
equation at the example of the NLS soliton,
\begin{equation}
    u(x,t) = 
    \sqrt{a}\mbox{ sech}\sqrt{a}(x-ct)\exp\left[i\left(\frac{c}{2}x+\left(a-\frac{c^{2}}{4}\right)t\right)\right]
    \label{sol}
\end{equation}
for $a=2$ and $c=15$ and $t\leq 2$ shown in Fig.~\ref{nlssol}. The exact solution (\ref{sol}) is 
imposed as initial condition at $t=0$, and the numerical error 
$\Delta$ is defined as in (\ref{Delta}). We use $x_{r}=-x_{l}=25$ to 
assure that the solution is almost of the order of machine precision at the 
boundaries (here $\sim 10^{-11}$). This is important since both the 
TBC and the PML approach in the form used here assume initial data 
with  compact support (at least within 
numerical precision) in the domain $[x_{l},x_{r}]$. As in 
\cite{Zhengpml}, the PML approach for NLS does not dissipate the 
solution well in the layer which leads to errors in the computational 
zone. 

As can be seen in Fig.~\ref{nlssol}, 
the localized solution travels with constant speed $c=15$ to the 
right. Its maximum  hits the boundary of 
domain II at $t=1.666\ldots$. The real part in the same 
figure shows the oscillatory phase of the solution that is stationary 
in a comoving frame (obtained by applying (\ref{eq:galilei})). In 
such a frame, Fourier methods with time integrators as in \cite{etna} 
would be the ideal choice, but as in the previous section, the goal 
is here to study the efficiency of the respective method when the 
soliton crosses the boundary. 

\begin{figure}[htb!]
   \includegraphics[width=0.49\textwidth]{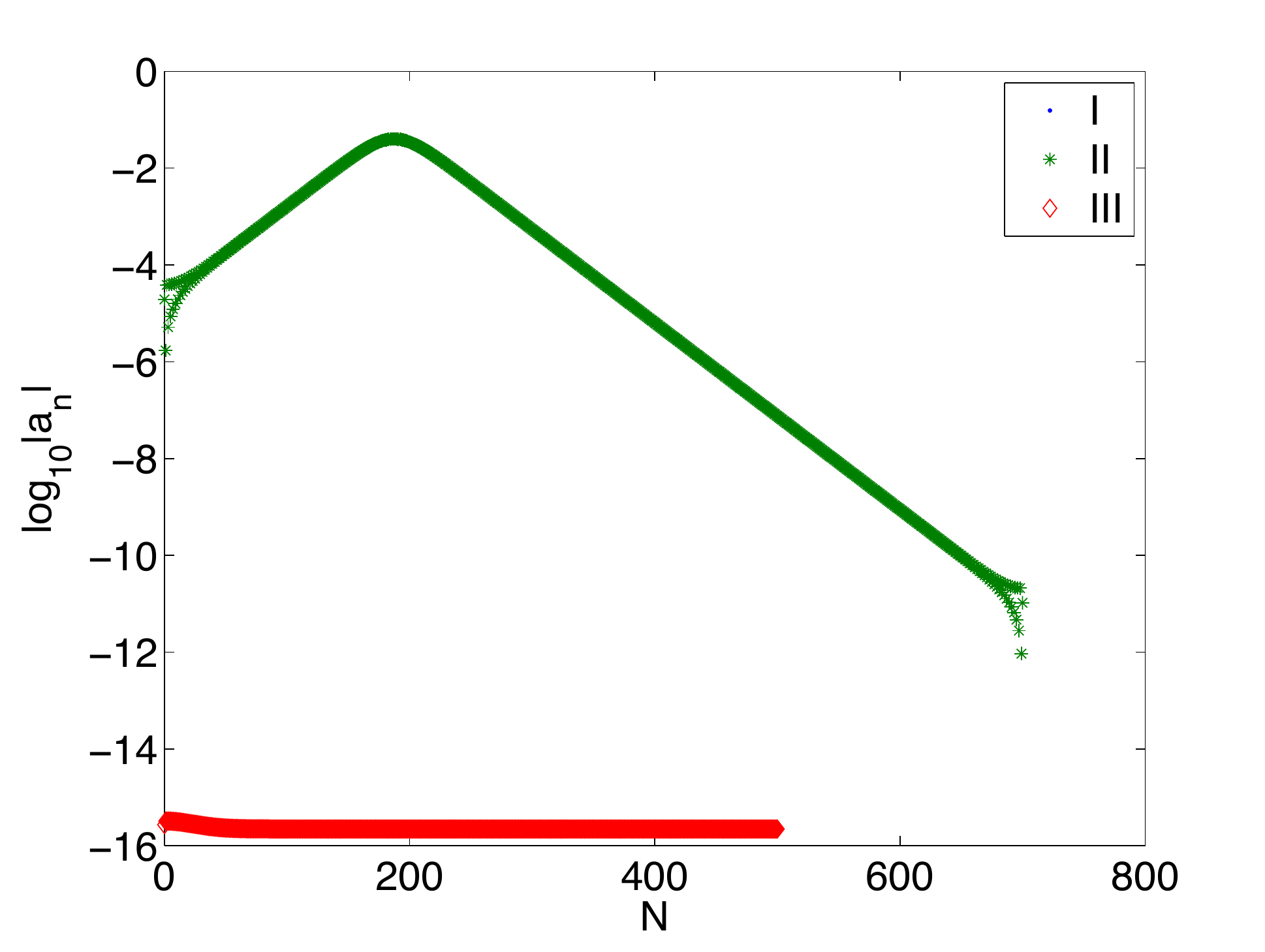}
   \includegraphics[width=0.49\textwidth]{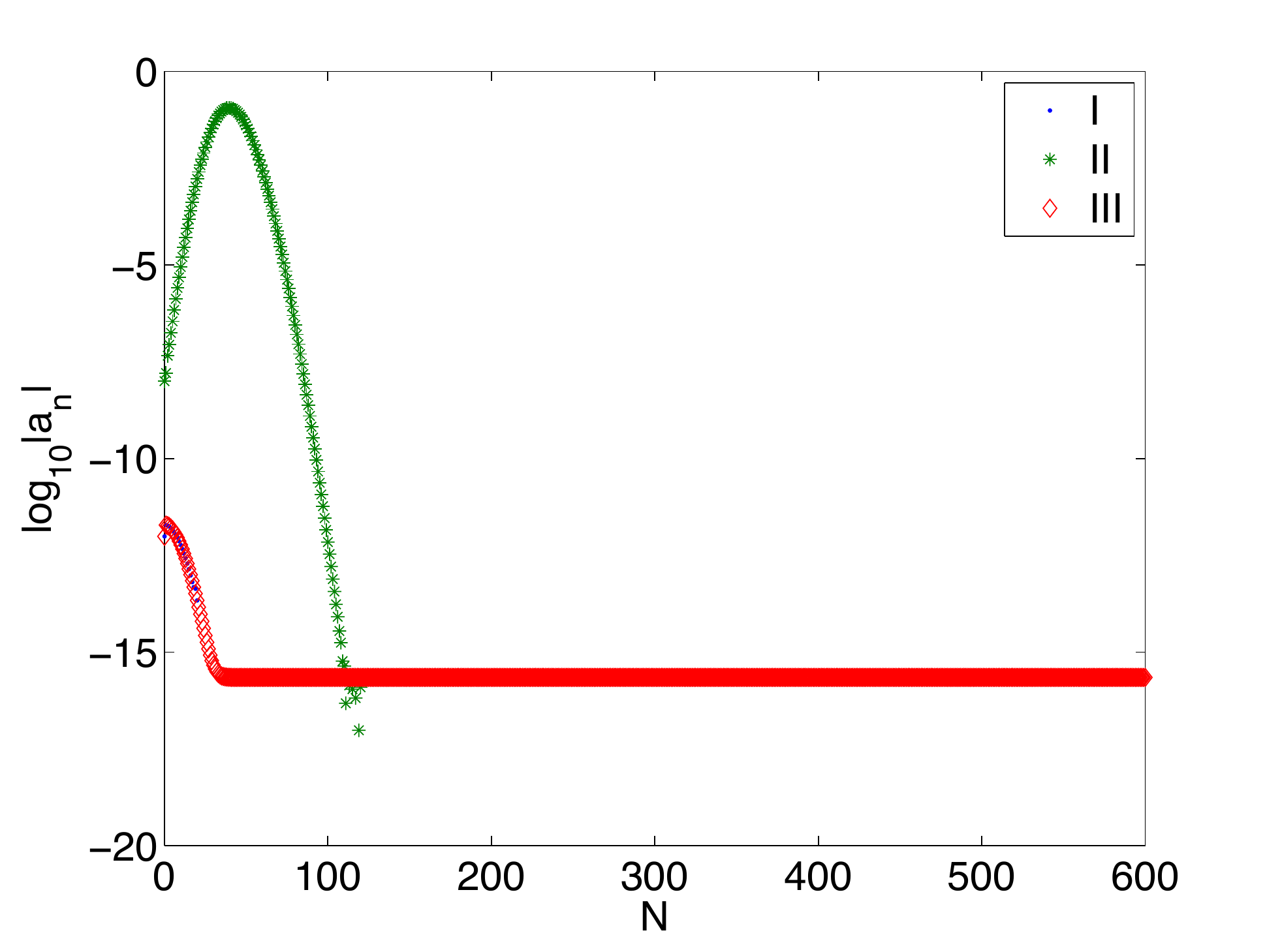}
 \caption{Chebyshev coefficients (\ref{coll}) for the  solution 
 (\ref{sol}) in the three domains; on the left at  $t=0$, on the 
 right for $t=2$.}
 \label{nlsres}
\end{figure}

We again study the spatial resolution for the CED method and choose 
$N^{I}=20$, $N^{II}=700$ and 
$N^{III}=500$ Chebyshev points in the respective domains. It can be 
seen in Fig.~\ref{nlsres} that the Chebyshev coefficients  
(\ref{coll}) decrease in each of the domains at least to $10^{-12}$ 
for this choice, 
both for the initial and the final time of the computation. Note that 
the Chebyshev coefficients in domain I are again completely of the order of 
machine precision  which means 
one could just impose vanishing Dirichlet condition at $x_{l}$. 

\begin{figure}[htb!]
   \includegraphics[width=0.49\textwidth]{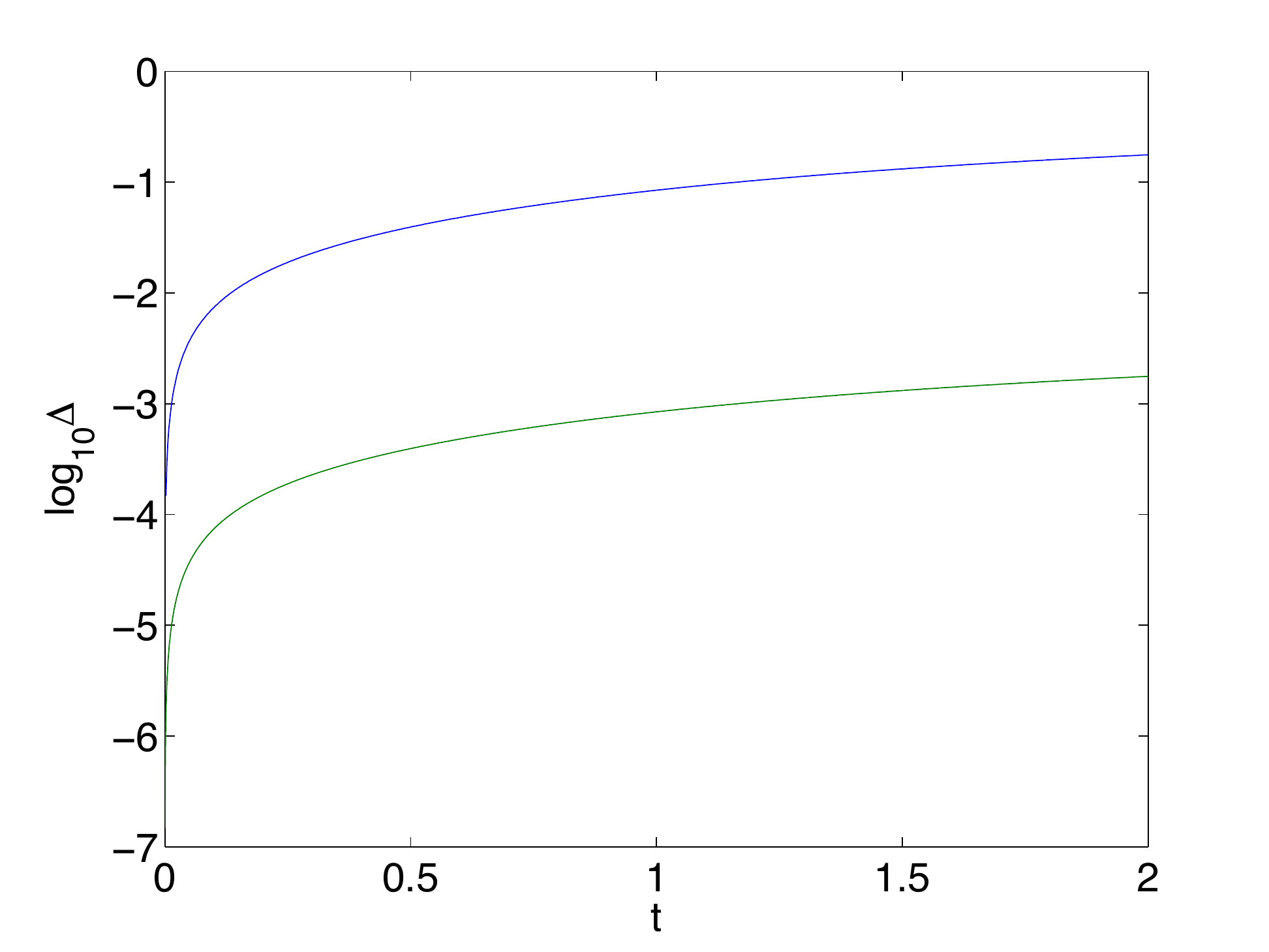}
   \includegraphics[width=0.49\textwidth]{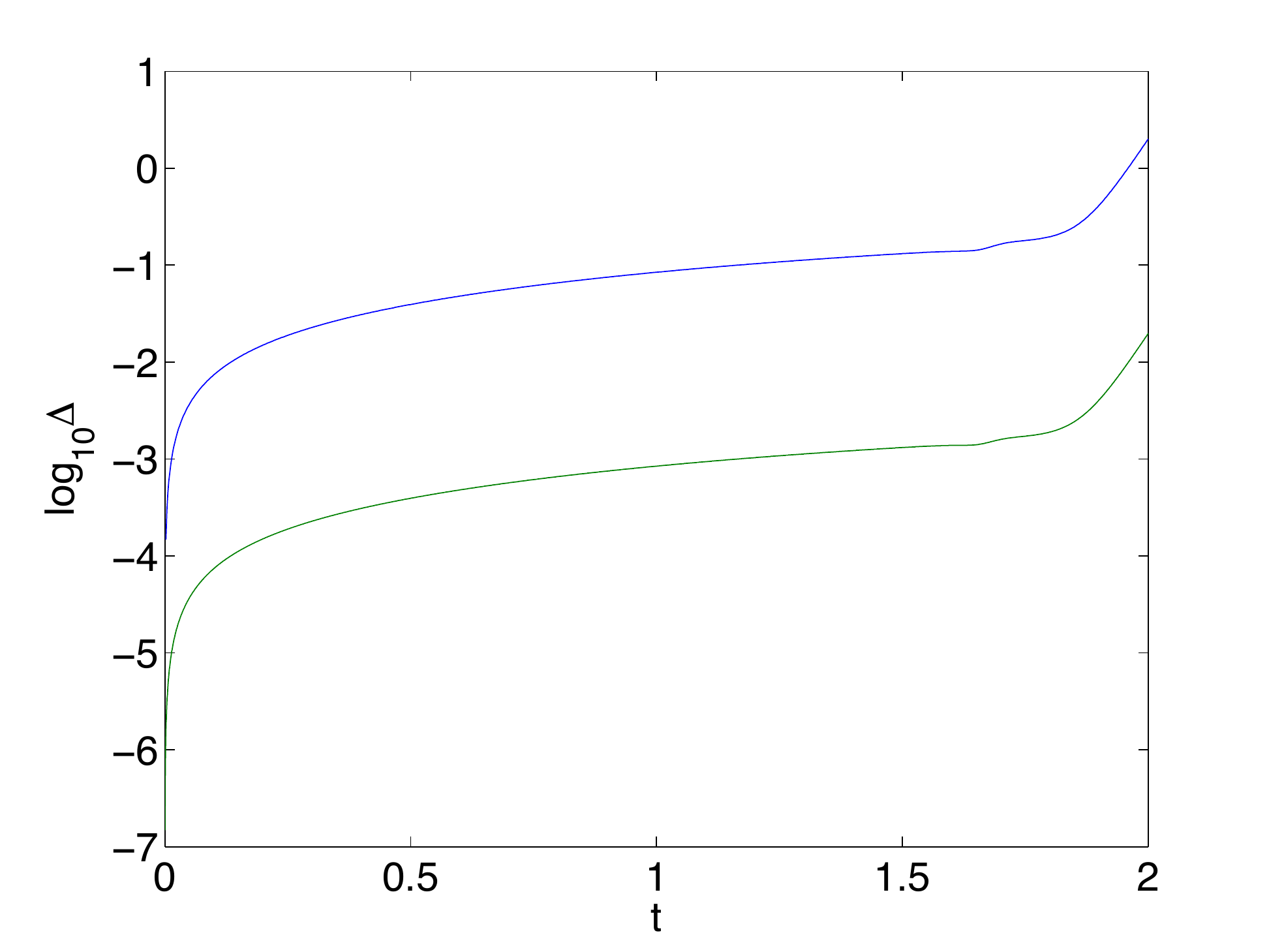}
 \caption{Normalized $L^{2}$ norm $\Delta$ (\ref{Delta}) of the difference between numerical 
 and exact solution with CN method for $N_{t}=10^{3}$ 
 and  
 $N_{t}=10^{4}$  time steps (from top to bottom); 
 on the left for the CED method the error on the whole line, on the 
 right for the TBC method in the computational domain.}
 \label{nlsmd}
\end{figure}

As in the previous section, we 
first compare the various methods with the second order CN scheme in 
all cases. The numerical error $\Delta$ for the CED method can be seen 
in Fig.~\ref{nlsmd} on the left. The global character of the approach 
ensures once more that  no effect is visible in the error
when the maximum crosses the boundary at $t=1.666\ldots$. Since all 
approaches are now iterative due to the nonlinearity of the NLS 
equation  in contrast to the linear case of the previous section,  
we only consider the resolutions 
$N_{t}=10^{3}$ and $N_{t}=10^{4}$. It can be clearly seen 
that the CN method shows the expected second order behavior. A 
precision of the order of $10^{-3}$ can be reached in this example 
without problems, 
for higher accuracies a fourth order method as discussed below would 
be clearly the better choice. 

In Fig.~\ref{nlsmd}, the same situation is studied on the right 
for the TBC approach. There is only domain II in this case where we use the 
same resolution as for CED ($N^{II}=700$). Due to the much lower 
overall 
spatial resolution needed than in the CED case, the method is slightly more 
efficient. But this is partly offset by the need to iterate the 
Dirichlet to von Neumann map. Whereas this has little effect for 
$x=x_{l}$ since the solution vanishes with numerical precision there, 
the iteration takes more and more time once the maximum of the 
solution approaches the boundary $x=x_{r}$. In addition the numerical 
error normalized by the $L^{2}$ norm of the exact solution within the 
computational domain now increases strongly once the maximum of the 
solution passes the boundary. For $N_{t}=10^{3}$ it is even of order 
1, for $N_{t}=10^{4}$ of the order of $1\%$.  The error normalized by 
the $L^{2}$ norm of the initial data still decreases to the order of 
$10^{-3}$ and $10^{-5}$ respectively. But the increase of the error 
in Fig.~\ref{nlsmd} is in clear contrast to the linear case in 
Fig.~\ref{schroemd} where the error even decreases. The reason for 
this appears to be that the Dirichlet to von Neumann map has to be 
obtained also approximatively, whereas it is explicitly known in the 
linear case. The CN scheme though
shows the expected second order behavior in both cases, it is just 
computationally expensive to 
reach higher precision with a second order scheme.

\begin{figure}[htb!]
   \includegraphics[width=0.49\textwidth]{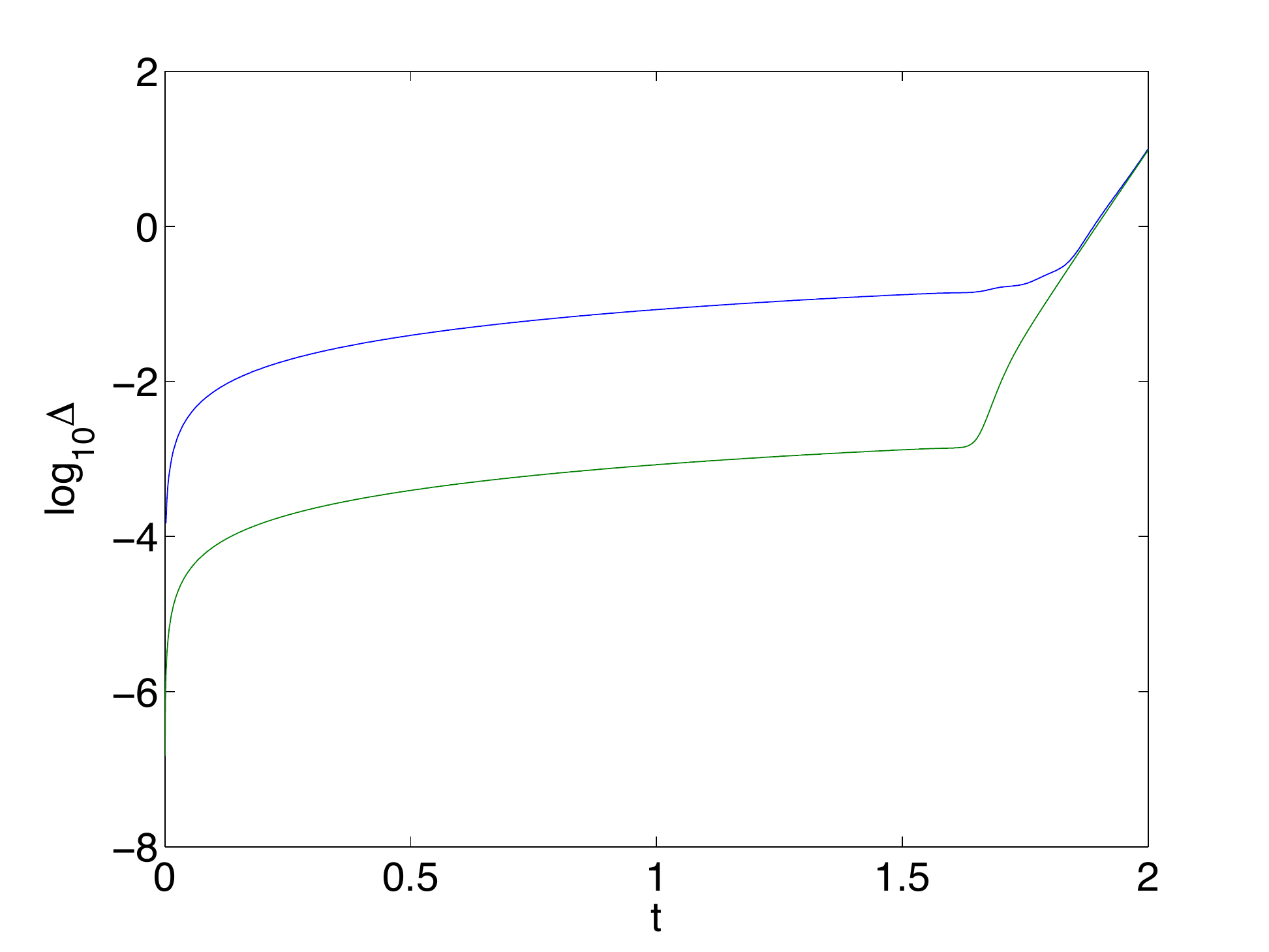}
   \includegraphics[width=0.49\textwidth]{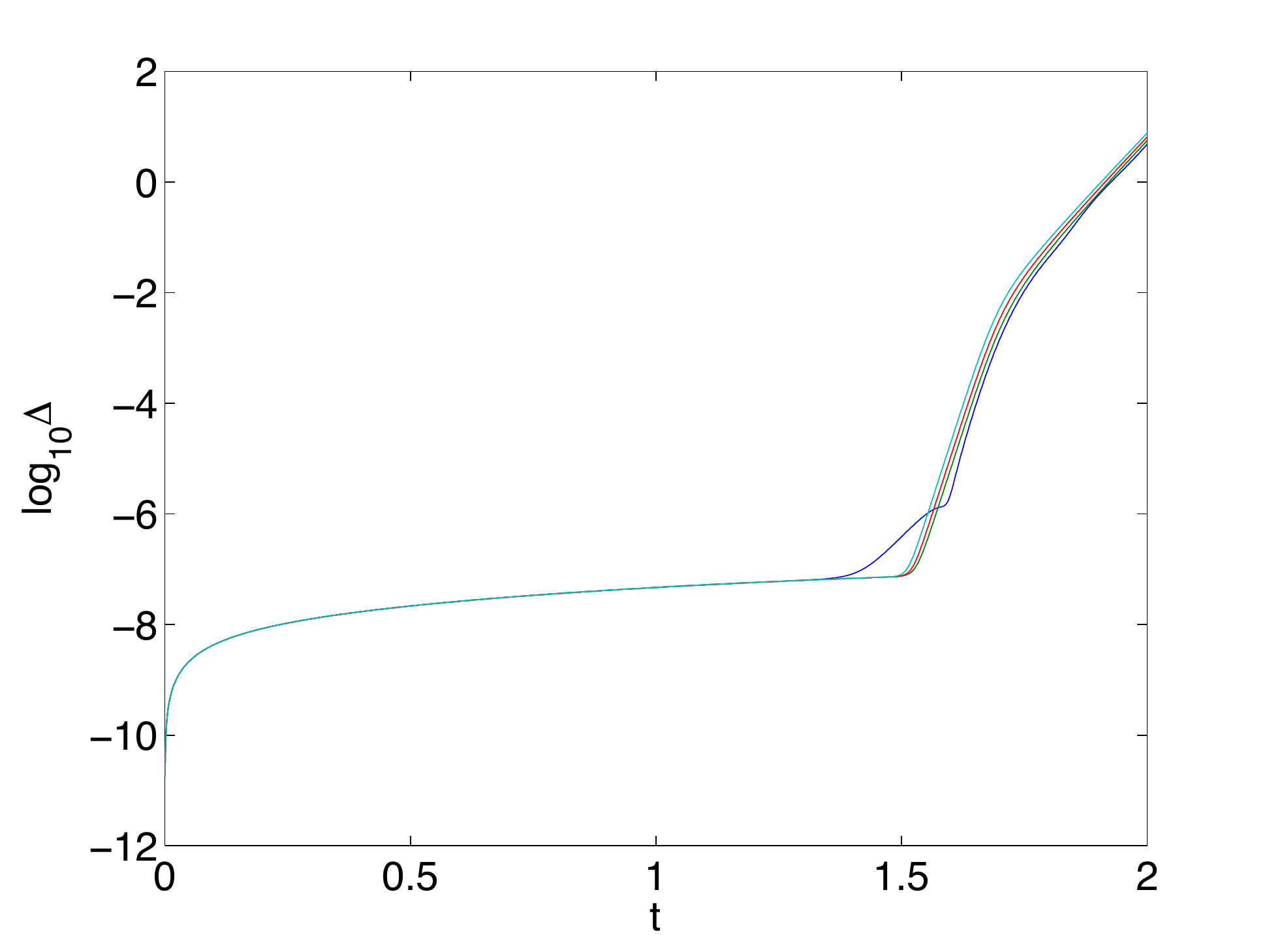}
 \caption{On the left, normalized $L^{2}$ norm $\Delta$ of the difference between numerical 
 and exact solution (\ref{Delta}) with CN method for $N_{t}=10^{3}$ 
 and 
 $N_{t}=10^{4}$ time steps (from top to bottom) 
 for the PML method with $\sigma_{0}=3$ and $\delta=1$; on the right the PML approach with IRK4 method 
 and $N_{t}=5000$ for various values of $\sigma_{0}$ in 
 (\ref{sigma}) ($\sigma_{0}=2,3,5,10$ corresponding to the blue, 
 green, red and cyan curve respectively).}
 \label{nlspml}
\end{figure}

As can be seen in Fig.~\ref{nlspml} on the left, the numerical 
error in the PML case with CN is obviously different from the 
TBC case. Again second order convergence of the time integration 
scheme is observed whilst most of the mass is within the 
computational domain, but the numerical error (\ref{Delta}) increases 
once most of the mass leaves the latter, and this even in a way 
independent of the time resolution at a given point. We use here 
$N^{I}=50$ and $N^{III}=100$ Chebyshev collocation points to assure 
that the Chebyshev coefficients (\ref{coll}) decrease to 
machine precision throughout the computation in the layers. Again a
much lower values of $N^{III}$ is needed here 
compared to the CED approach. 

\begin{figure}[htb!]
   \includegraphics[width=0.49\textwidth]{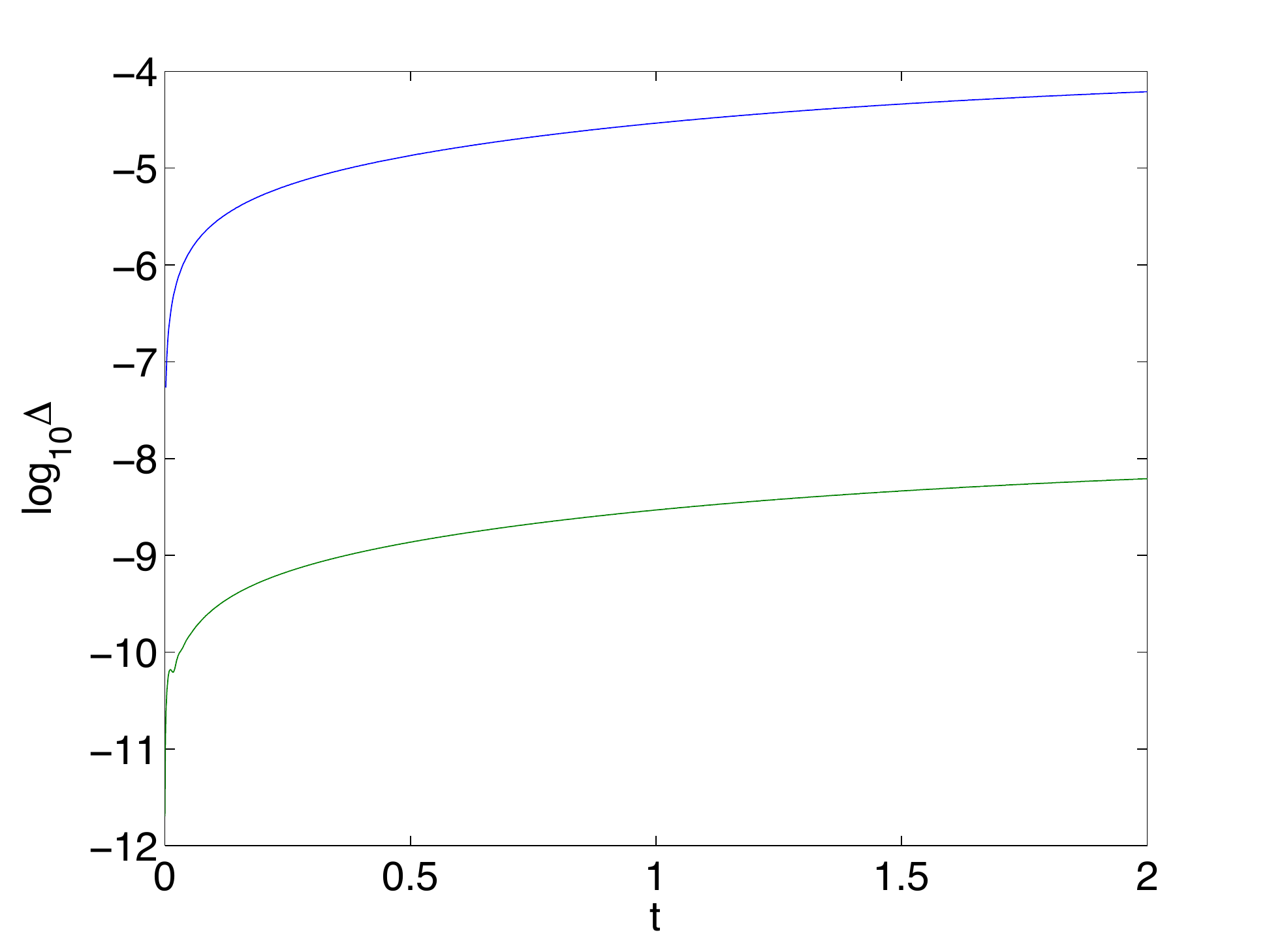}
   \includegraphics[width=0.49\textwidth]{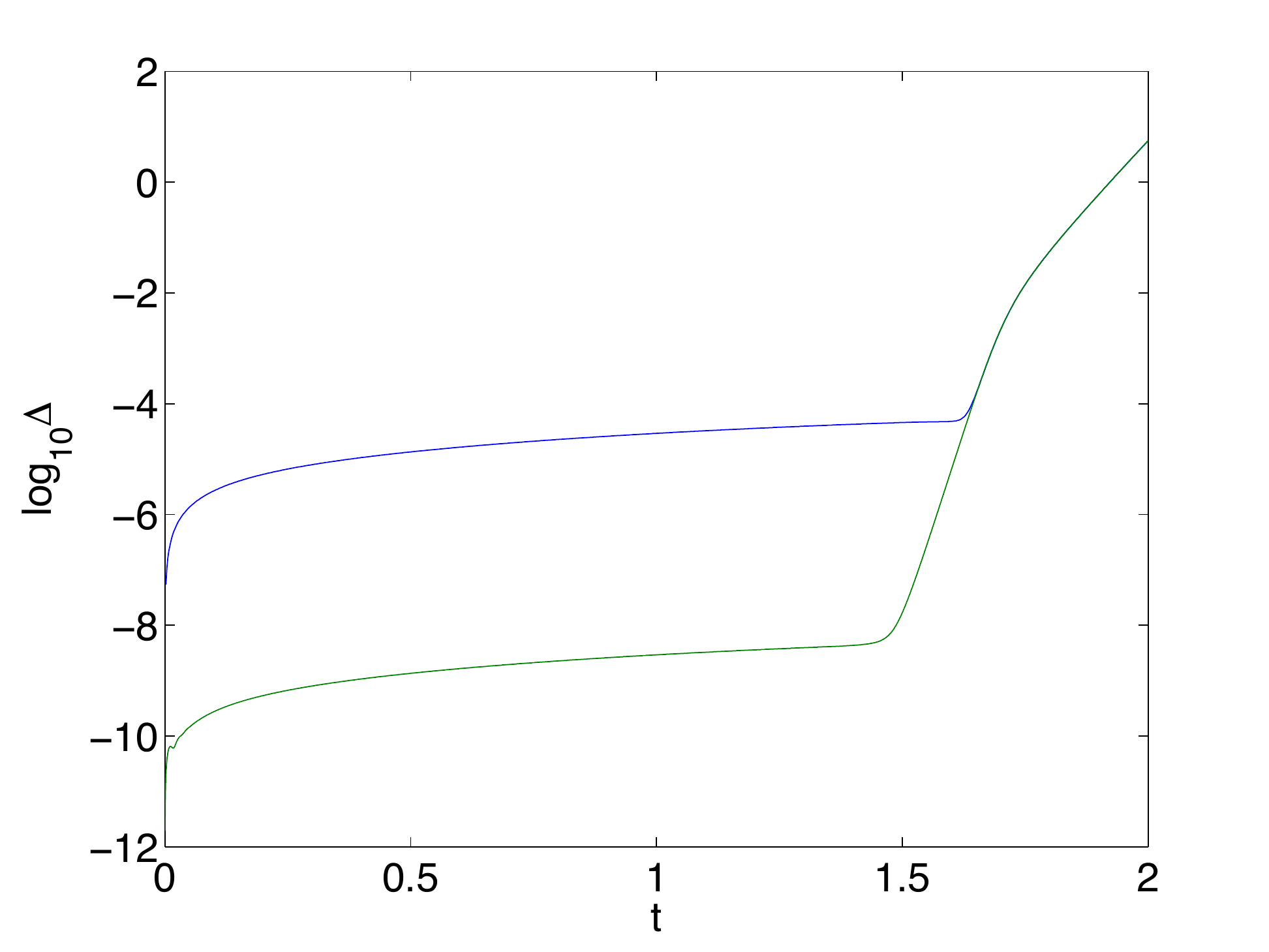}
 \caption{Normalized $L^{2}$ norm $\Delta$ of the difference between numerical 
 and exact solution (\ref{Delta}) with IRK4 method for $N_{t}=10^{3}$ 
 and  
 $N_{t}=10^{4}$ time steps (from top to bottom); 
 on the left for the CED method the error on the whole line, on the 
 right for the PML method with $\delta=1$ and $\sigma_{0}=3$ in the computational domain.}
 \label{nlsirk4}
\end{figure}

To reach higher precision we apply again 
the fourth order IRK4 method which can be implemented for both PML and 
CED in a straight forward way. The fourth order convergence 
can be seen in Fig.~\ref{nlsirk4}. For CED, with $N_{t}=10^{4}$ one is 
already close to  machine precision. 

For the PML method, the same situation is shown on the right of 
Fig.~\ref{nlsirk4}. As with the CN method in Fig.~\ref{nlspml}, the 
expected convergence of the scheme is only observed whilst most of 
the mass of the solution is within the computational domain. Once the 
maximum of the solution is close to the boundary, the error becomes 
independent both of the spatial and the temporal resolution and is 
thus clearly an effect of the PML approach. As observed already in 
\cite{Zhengpml}, the layers are clearly not `perfectly matched'  in the nonlinear case. 
This also does not change if we vary the parameter $\sigma_{0}$ in 
the right figure of Fig.~\ref{nlspml}. There 
we put $\delta=1$ and vary $\sigma_{0}$ to minimize the 
numerical error. It can be seen in Fig.~\ref{nlspml} that the 
numerical error $\Delta$ is almost identical for the values 
$2,3,5,10$, but that it is worst for the values 2 and 10 in this 
example. Thus we choose $\sigma_{0}=3$ since it also performs best 
close to the time when the maximum of the solution reaches the 
boundary of the computational domain. 

\begin{figure}[htb!]
   \includegraphics[width=0.49\textwidth]{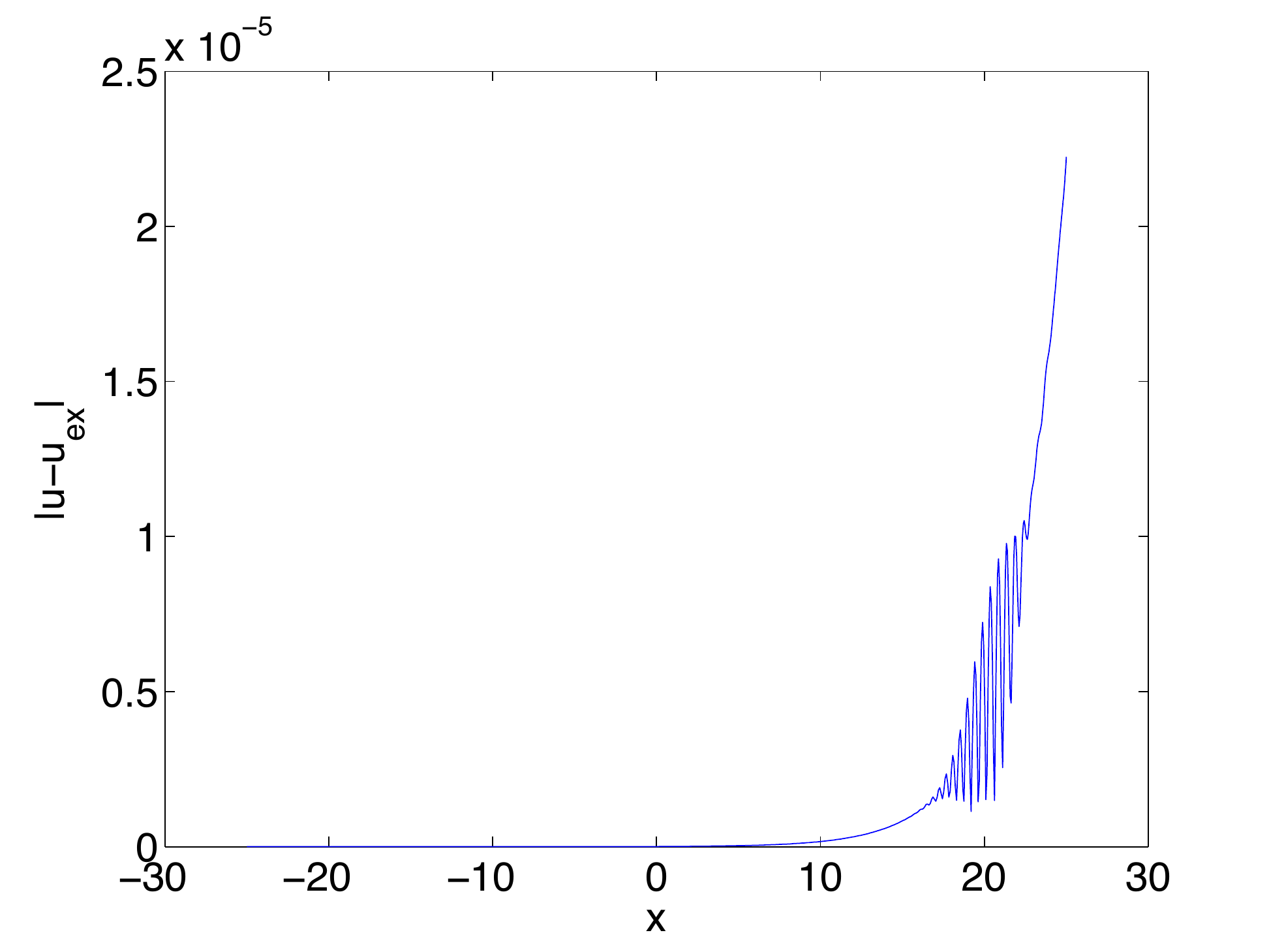}
   \includegraphics[width=0.49\textwidth]{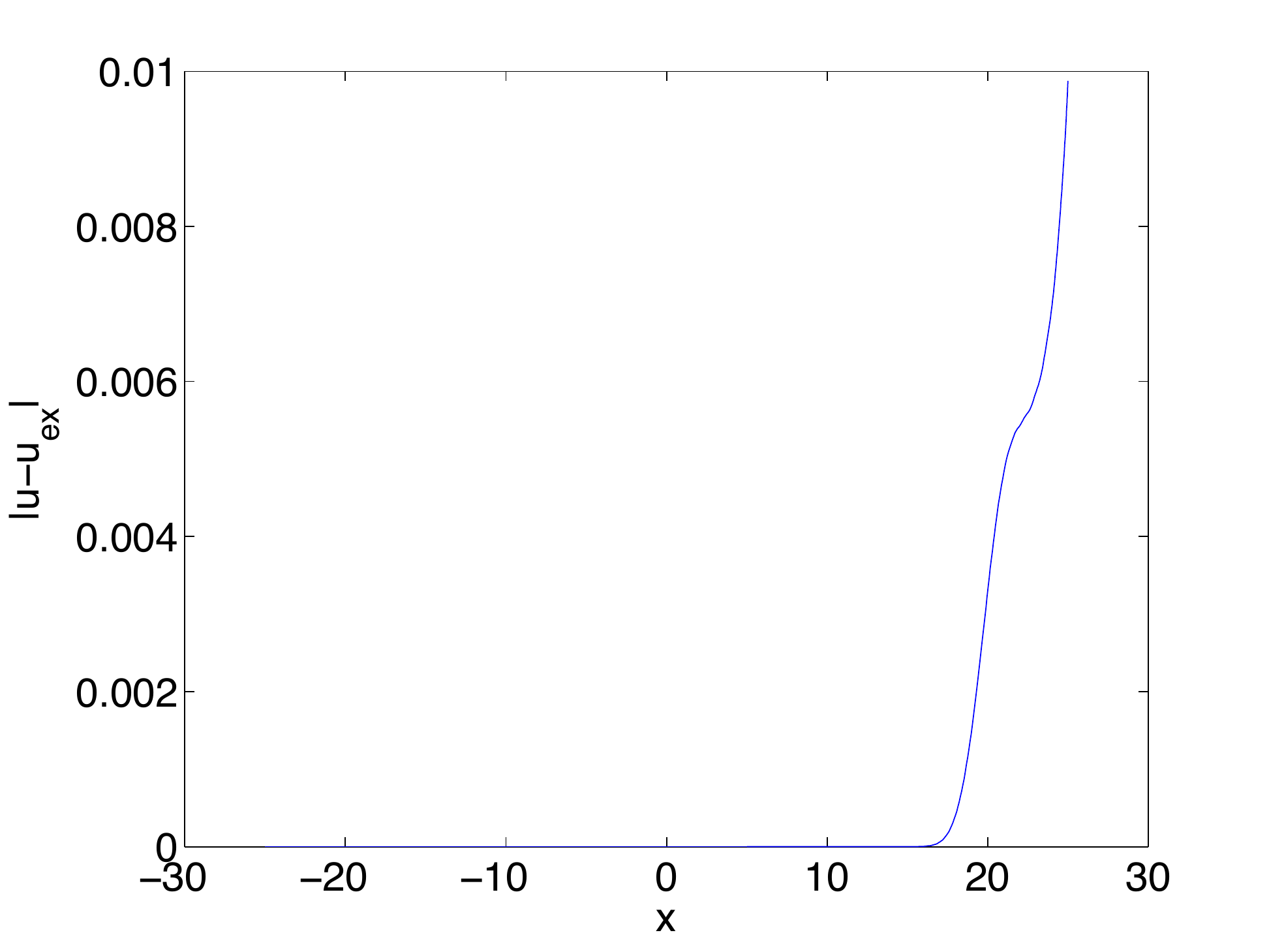}
 \caption{Difference between numerical 
 and exact solution to NLS for soliton initial data (\ref{sol}) for the TBC method on the left and 
 for the PML method on the right at $t=2$ for 
 $N_{t}=10^{4}$ with the CN method.}
 \label{nlsdiff}
\end{figure}

To illustrate the behavior of 
the numerical solution for TBC and PML near the boundary, we show in 
Fig.~\ref{nlsdiff} the 
difference of the numerical and exact solutions in both cases at 
$t=2$ for CN with $N_{t}=10^{4}$. It can be seen that the absolute error for TBC is of the order 
$10^{-5}$ and that some spurious reflections propagate to the left. 
These might be further reduced if a higher resolution in time is 
used. The absolute error for PML on the other hand is of the order of 
$1\%$, and the solution in the computational domain is clearly 
affected near the boundary. This effect cannot be reduced if the 
resolution in time or space is increased, nor if different values of 
$\sigma_{0}$ are used. 

\section{Numerical study of the Peregrine breather}
\begin{figure}[htb!]
   \includegraphics[width=0.7\textwidth]{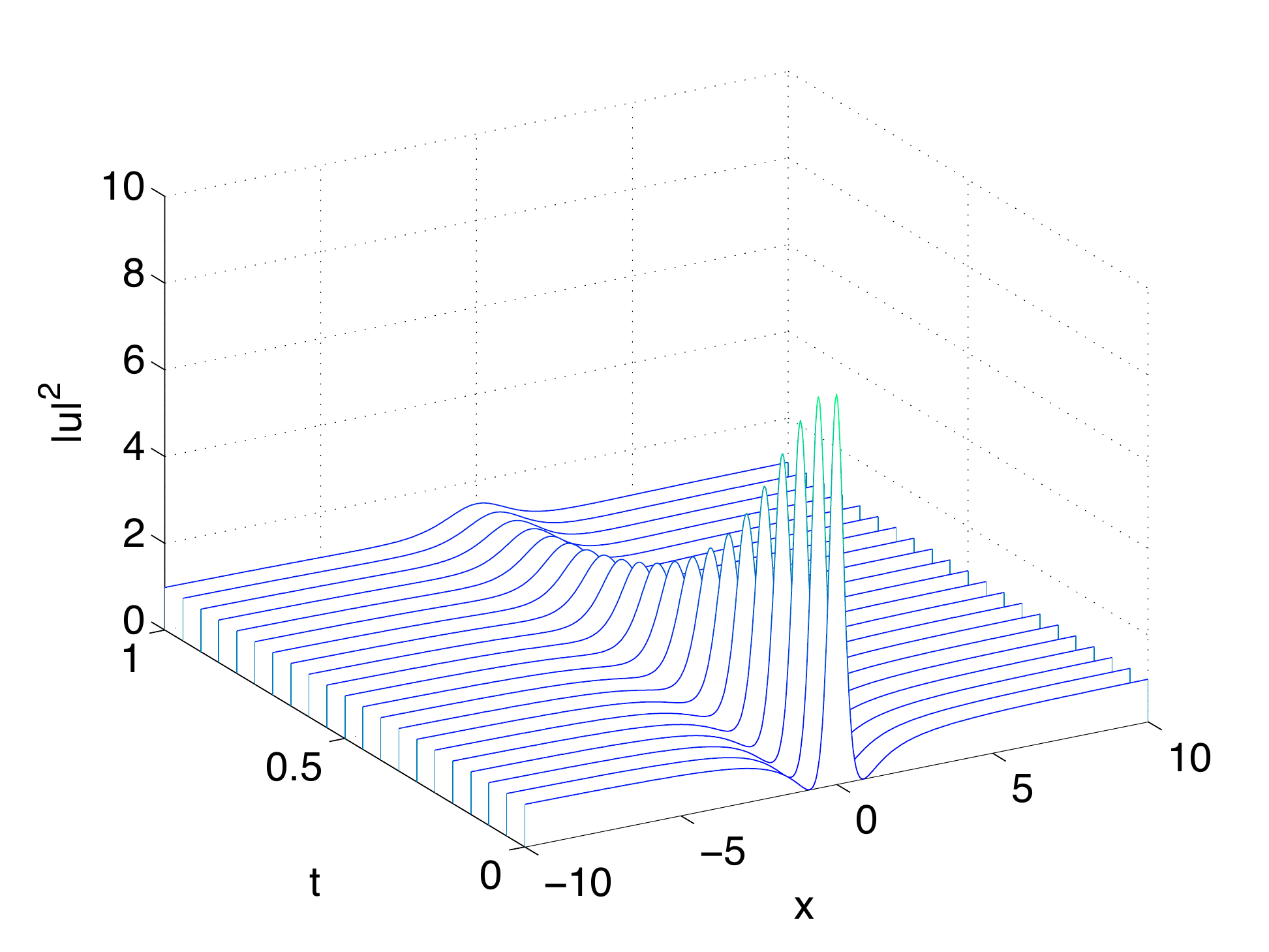}
 \caption{Peregrine breather (\ref{peregrine}).}
 \label{breatherfig}
\end{figure}

The TBC and the PML method used in this article are constructed for 
initial data with compact support in the computational domain, a 
condition which should be satisfied with the aimed at precision. CED 
methods on the contrary just require  fall off to some (possibly 
nonzero) constant (for fixed $t$) at infinity which may be slow. To illustrate this 
we consider the celebrated Peregrine breather solution to the NLS 
equation which is discussed as a possible candidate for rogue 
waves in hydrodynamics and nonlinear optics. 
The solution can be seen in Fig.~\ref{breatherfig}.
For $|x|\to\infty$, one has that $|u|\to1$. The solution is also 
slowly dispersed away for $t\to\infty$. We show in this section 
that the breather solution can be numerically evolved on the whole 
real line with essentially machine precision. This allows to study 
numerically the stability of the breather solution. This is important 
in the context of rogue waves since experimentally observable 
structures should be stable to a certain extent. A stability analysis 
for periodic breathers was presented in \cite{schober}, an 
experimental study of the stability of the Peregrine breather 
appearing in hydrodynamical experiments against wind was described in 
\cite{cha3}.

\begin{figure}[htb!]
   \includegraphics[width=0.49\textwidth]{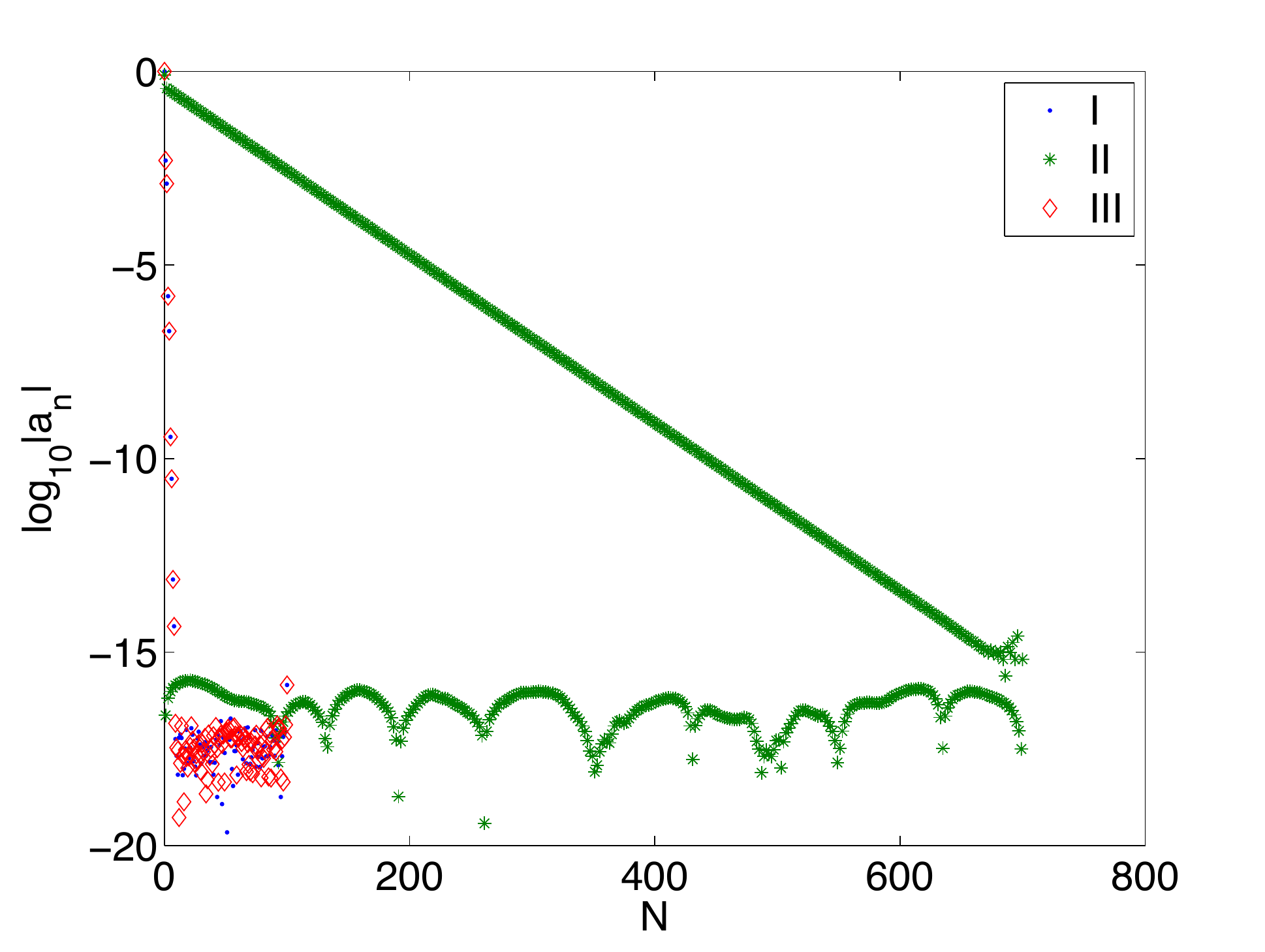}
   \includegraphics[width=0.49\textwidth]{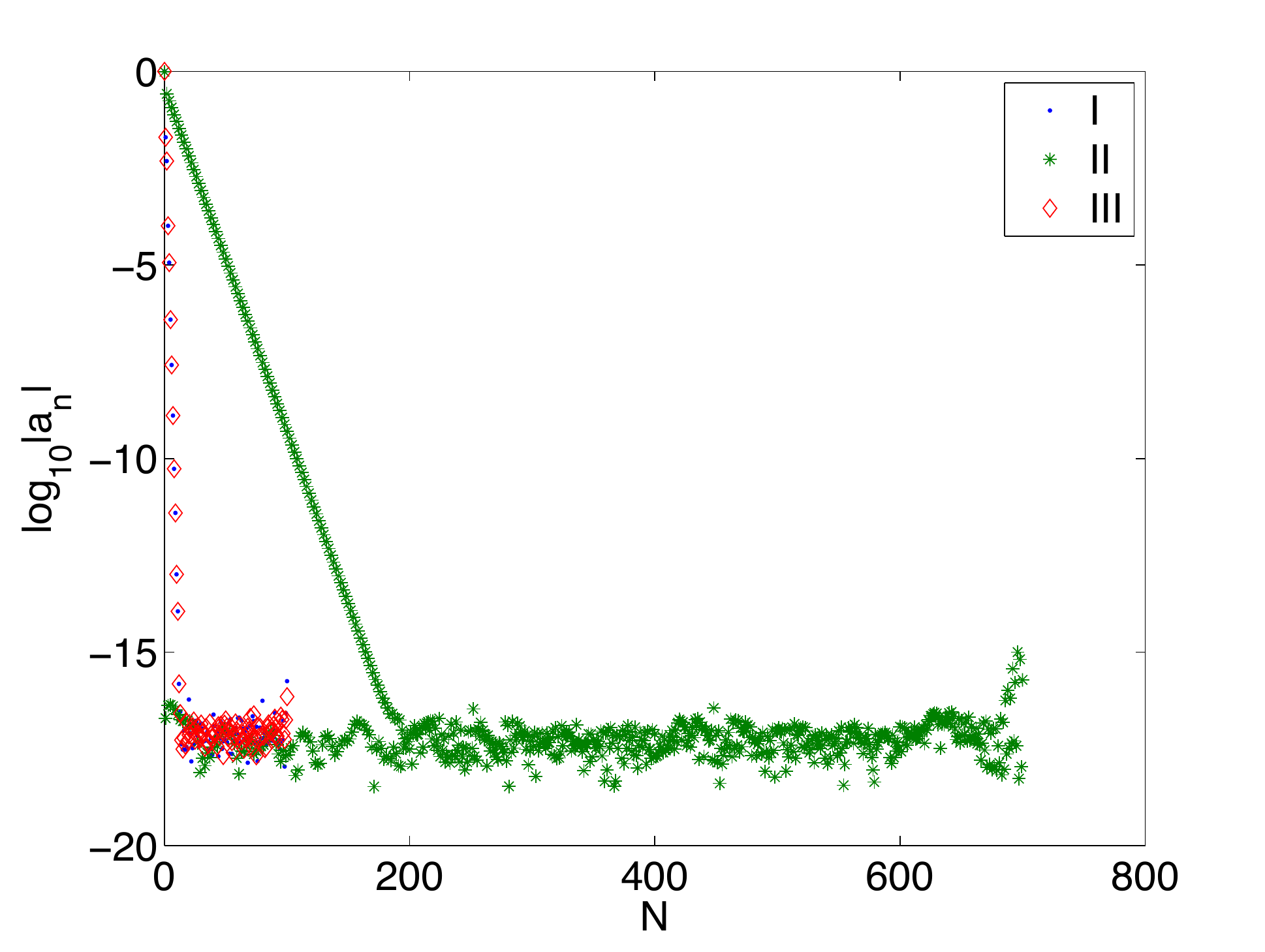}
 \caption{Chebyshev coefficients (\ref{coll}) for the  solution 
 (\ref{peregrine}) in the three domains; on the left at  $t=0$, on the 
 right for $t=1$.}
 \label{breathercheb}
\end{figure}

We give the Peregrine solution for $t=0$ as initial data and 
numerically solve NLS with these data for $t\leq 1$. To this end we 
choose $x_{r}=-x_{l}=10$, $N^{I}=N^{III}=50$ and $N^{II}=700$. It can 
be seen in Fig.\ref{breathercheb} that the Chebyshev coefficients (\ref{coll}) decrease to 
machine precision in this case during the whole computation. 

\begin{figure}[htb!]
   \includegraphics[width=0.7\textwidth]{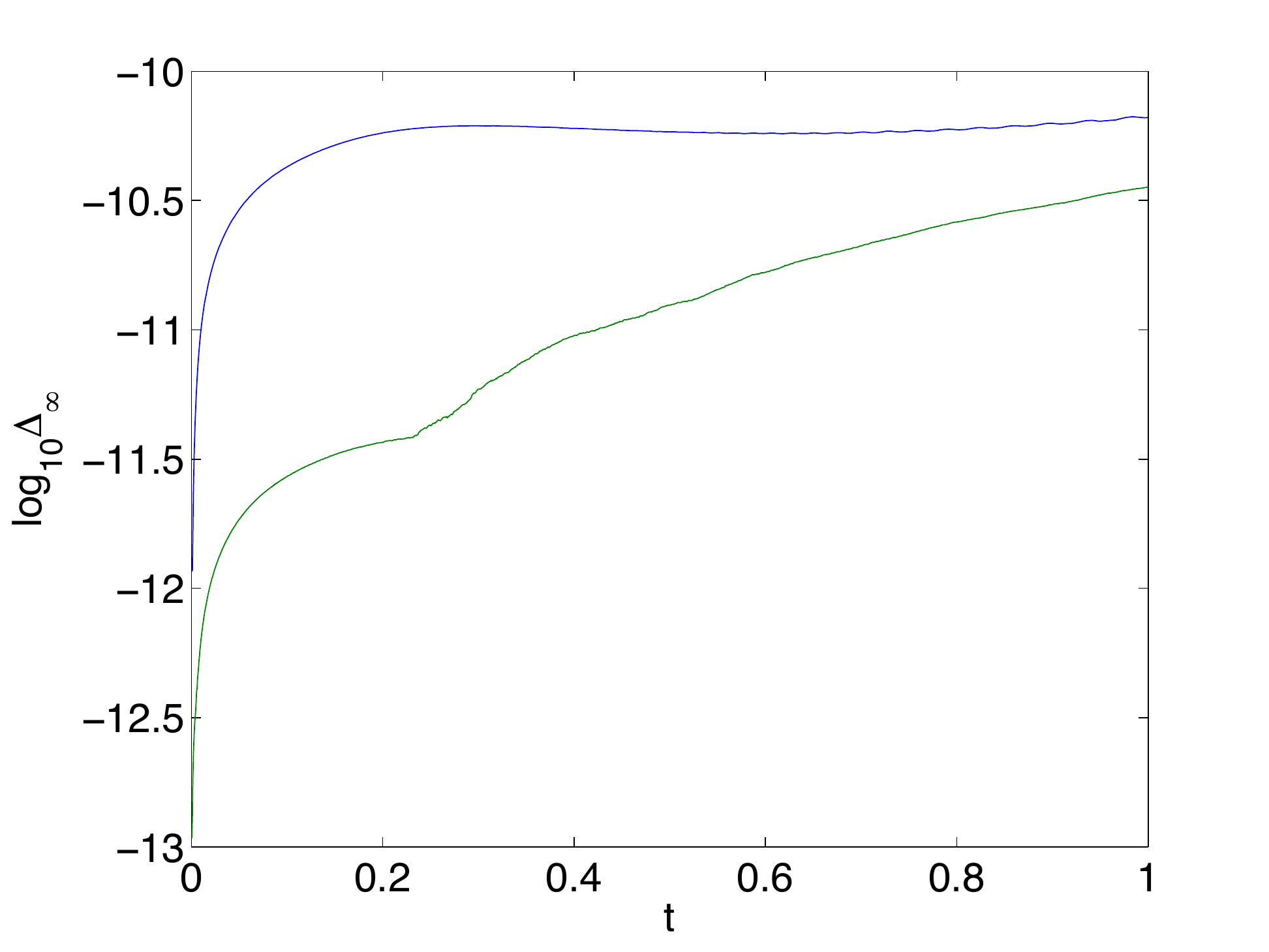}
 \caption{$L^{\infty}$ norm $\Delta_{\infty}$ (\ref{Deltainf}) of the 
 difference between numerical 
 and exact solution to NLS for Peregrine initial data 
 (\ref{peregrine}) for the CED approach with IRK4 method with 
 $N_{t}=1000$ and $N_{t}=2000$ (from top to bottom).}
 \label{breather}
\end{figure}

Note 
that we get essentially the same numerical result with the more 
symmetric (in the $N^{i}$) choice $x_{r}=-x_{l}=5$ and 
$N^{I}=N^{III}=200$, $N^{II}=400$. We use here only the IRK4 method 
since we are interested in testing with which accuracy the solution 
can be reproduced. Since the Peregrine solution is not in 
$L^{2}(\mathbb{R})$, we consider here the numerical error defined via 
the $L^{\infty}$ norm,
\begin{equation}
    \Delta_{\infty}=\frac{||u-u_{Per}||_{\infty}}{||u_{Per}||_{\infty}},
    \label{Deltainf}
\end{equation}
where $u_{Per}$ is the solution (\ref{peregrine}). 
It can be seen in Fig.~\ref{breather} that the numerical error is 
smaller than $10^{-10}$ for $N_{t}=1000$, and that one is close to the optimum with 
$N_{t}=2000$. 

If no exact solution is known for given initial data, it is 
convenient to use conserved quantities of the studied PDE, whose 
conservation is not implemented in the code, to control the accuracy 
of the numerical solution, see for instance the discussion in 
\cite{etna}. Since the cubic NLS equation is completely integrable, 
there is an infinite number of conserved quantities, the most popular 
being the mass (the $L^{2}$ norm of the solution) and the energy. 
Note that both these quantities are not defined for the breather 
since it is not in $L^{2}$, but that the combination 
\begin{equation}
    E=\frac{1}{2}\int_{\mathbb{R}}^{}\left\{|u_{x}|^{2}-|u|^{2}(|u^{2}|-1)\right\}dx
    \label{EB}
\end{equation}
is both conserved and defined. It can be 
checked by direct computation that $E=0$ for the Peregrine breather. For given 
$u$, the integral in (\ref{EB}) will be computed as discussed in 
section \ref{spectral} with the Clenshaw-Curtis algorithm after 
division by $s_{l}$ and $s_{r}$ in coefficient space in the 
compactified domains. In 
numerical simulations, the quantity $E$ will depend on time due to 
unavoidable numerical errors. We use the 
relative quantity $\Delta_{E}=|1-E(t)/E(0)|$ as an indicator of the 
accuracy of the solution. As discussed in \cite{etna}, such 
quantities are reliable indicators of the precision if sufficient 
resolution in $x$ is provided, but tend to overestimate the accuracy 
by one to two orders of magnitude. For the exact Peregrine solution, 
the quantity $E$ is numerically of the order of $10^{-13}$, i.e., 
machine precision.

\begin{figure}[htb!]
   \includegraphics[width=0.49\textwidth]{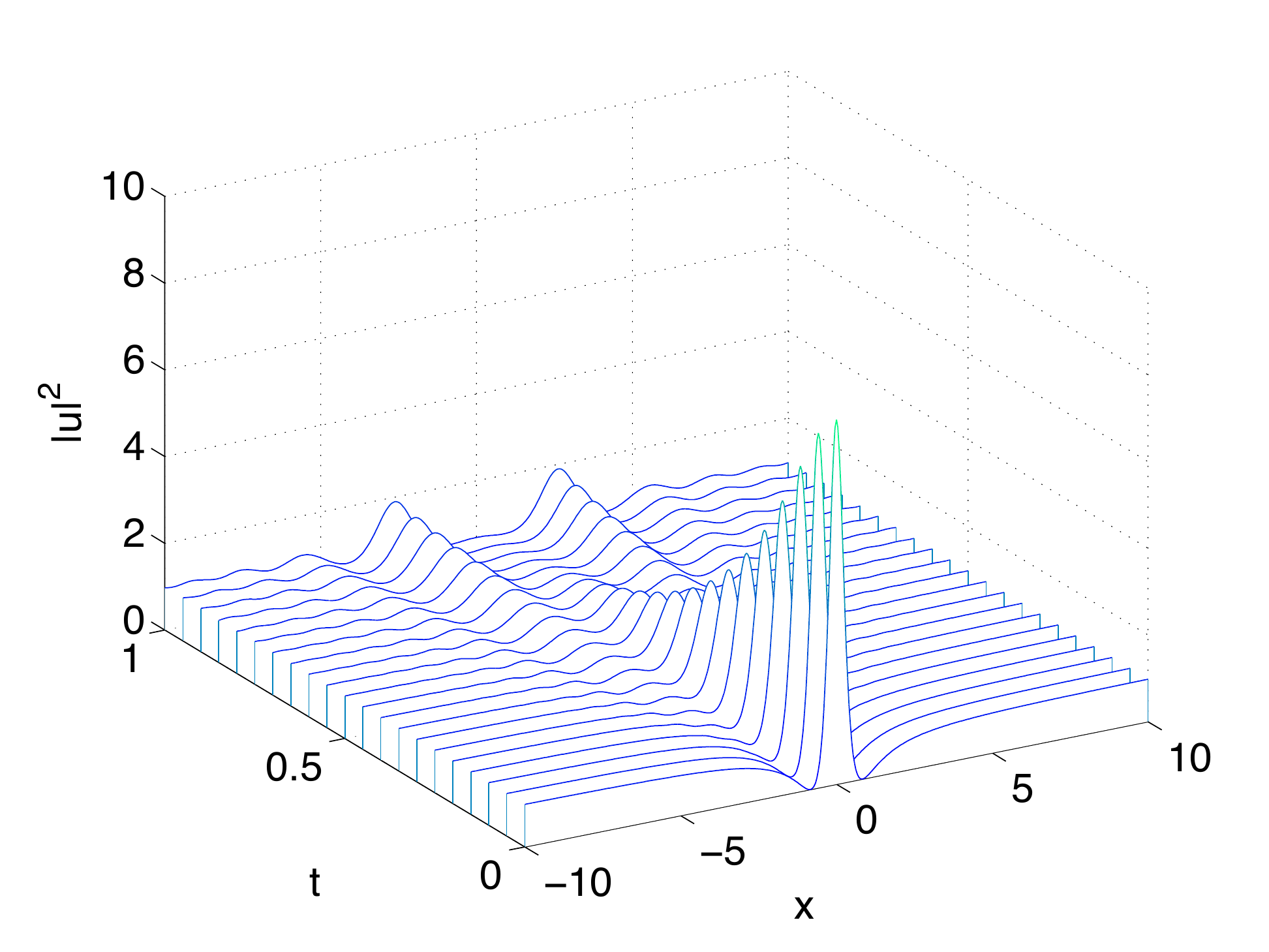}
   \includegraphics[width=0.49\textwidth]{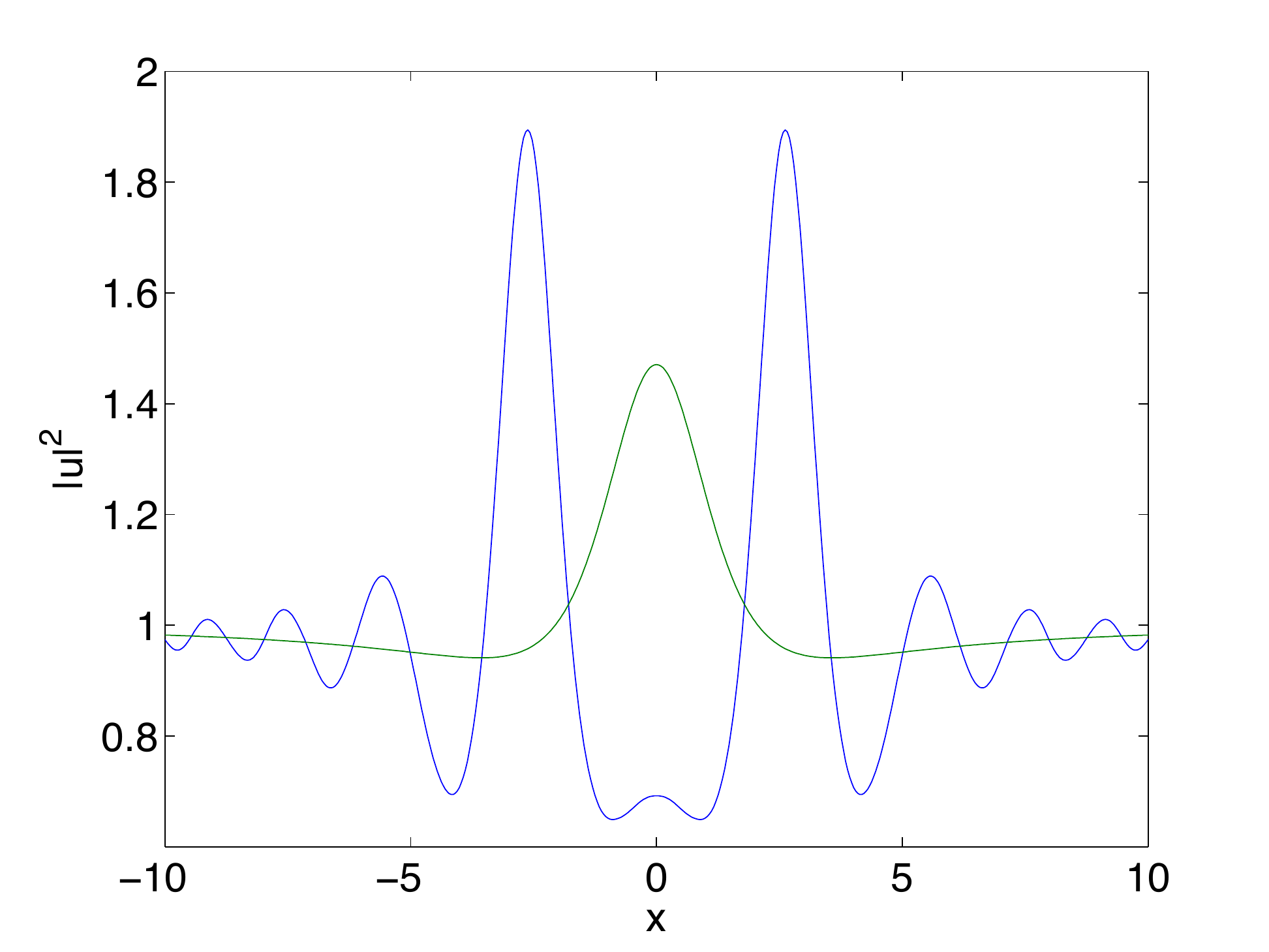}
 \caption{Solution to the NLS equation for the initial data 
 $u_{Per}(x,0)+0.1\exp(-x^{2})$ in dependence of time on the left; on 
 the right the solution at the final time in blue and the Peregrine 
 solution at the same time in green.}
 \label{breathergauss}
\end{figure}

Since the CED approach allows to propagate Peregrine initial data 
essentially with machine precision, we are able to study localized
perturbations of the breather: as an example we consider initial data 
$u_{Per}(x,0)+0.1\exp(-x^{2})$, i.e., the Peregrine solution at time 
$t=0$ perturbed by a small Gaussian. We use 
$N^{I}=N^{II}=N^{III}=400$ and $N_{t}=1000$ for 
the computation. It can be 
seen in Fig.~\ref{breathergauss} that the solution appears to be 
dispersed away to infinity, but in a significantly different way than 
the Peregrine breather. This is even more obvious in the right figure of 
Fig.~\ref{breathergauss} where the solution at the last computed time 
is shown with the Peregrine solution (\ref{peregrine}) at the same 
time.

\begin{figure}[htb!]
   \includegraphics[width=0.49\textwidth]{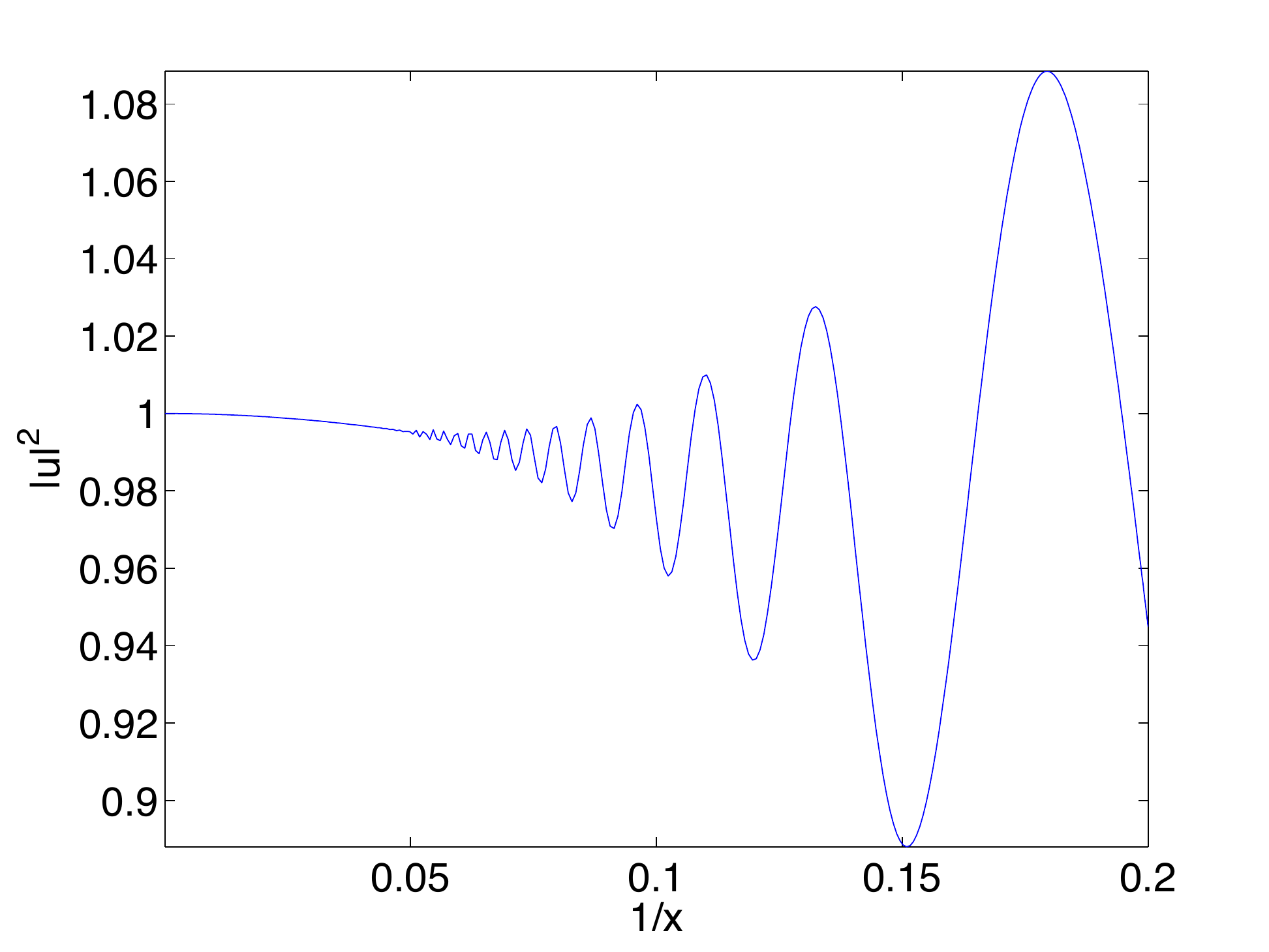}
   \includegraphics[width=0.49\textwidth]{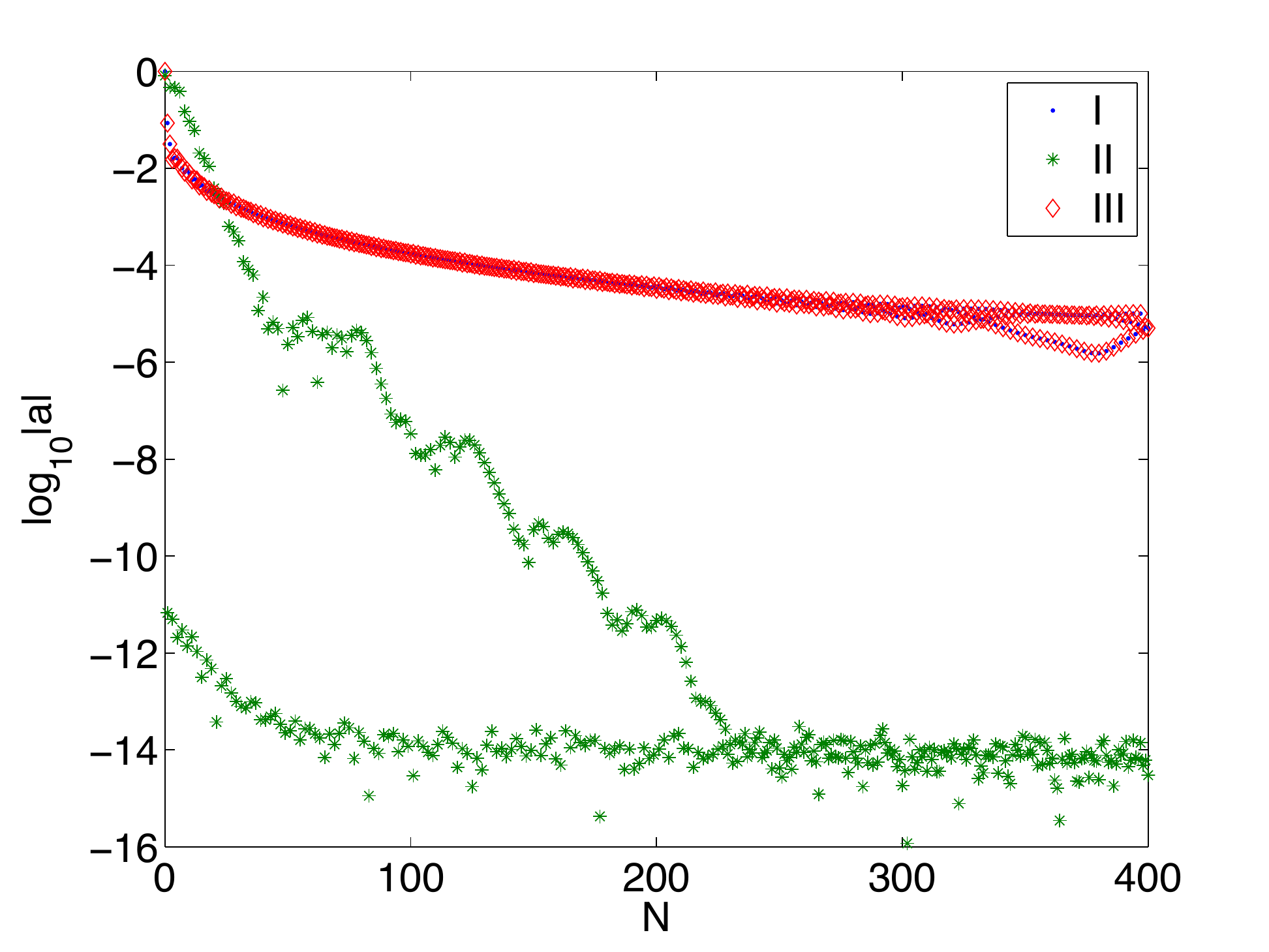}
 \caption{Solution to the NLS equation at $t=1$ for the initial data 
 $u_{Per}(x,0)+0.1\exp(-x^{2})$  in the compactified domain on the left, and the 
 corresponding Chebyshev coefficients (\ref{coll}) on the right.}
 \label{breathergausscheb}
\end{figure}

Thus it appears that the Peregrine breather is unstable against 
this type of perturbations, and that the perturbed solution does not 
stay close to the exact solution. Note that the solution is even at 
the last computed time well resolved spatially as can be seen in the 
right figure of Fig.~\ref{breathergausscheb}, where the Chebyshev 
coefficients (\ref{coll}) at the last computed time are shown. Since 
the initial data is symmetric with respect to the transformation $x\to-x$ and since the 
Schr\"odinger equation preserves parity, the Chebyshev coefficients 
in zone I and III are identical, and half of the coefficients in zone 
II vanish with numerical precision. The coefficients in zone I and 
III still decrease to $10^{-4}$ which indicates that the solution is 
computed to better than plotting accuracy. The slow decrease of the 
Chebyshev coefficients is due to the oscillations shown on the left 
of Fig.~\ref{breathergausscheb}. The relative conserved energy 
$\Delta_{E}\sim 8.8*10^{-3}$. This implies together with the 
resolution in coefficient space shown in 
Fig.~\ref{breathergausscheb} that the solution is computed to at 
least plotting accuracy. 

\section{Outlook}
In this paper we have presented a multidomain spectral method with 
compactified exterior domains for Schr\"odinger equations and 
discussed examples for the linear Schr\"odinger and cubic NLS 
equations. The solutions have to be bounded at spatial infinity. 
The method was compared to exact TBC  and PML 
approaches. It was shown that it produces results of at least the 
same quality as TBC and PML in the linear case, and that machine 
precision can be reached with a fourth order time integrator. The 
price to pay for this is that the solution has to be resolved with 
the wanted precision on the whole real line. Thus if one is only 
interested in the solution in the computational domain, TBC and PML 
approaches might be more economic, though the latter approach in 
practice requires several runs to optimize the parameters. The 
slight  disadvantage of CED 
disappears, however, in the nonlinear case where even TBC 
methods need an approximate (and iterative) solution of the Dirichlet 
to von Neumann map, and where the PML method \cite{Zhengpml} 
produces considerable errors. In contrast, the CED method could reach 
the same precision as in the linear case and seems therefore the 
preferred choice 
for high precision studies. 

The real advantage of the CED method is, however, for situations with 
an algebraic fall off towards infinity whereas TBC and PML require 
initial data with compact support. As an example we studied 
the Peregrine breather (\ref{peregrine}) which is discussed as a 
model for rogue waves in hydrodynamics and nonlinear optics. As was 
shown, the CED approach can propagate this solution with machine 
precision and allows the 
study of perturbations which will be done in more detail elsewhere. 
The approach is also open to a generalization to higher dimensions 
where soliton solutions often have an algebraic fall off, see the 
\emph{lump} of Davey-Stewartson and Kadomtsev-Petviashvili equations
($2+1$ dimensional generalizations of NLS and Korteweg-de Vries 
equations respectively). This will be the subject of further research.


\begin{thebibliography}{99}
\bibitem{AntoineArnold}
X. Antoine, A. Arnold, C. Besse, M. Ehrhardt,  A. Sch\"adle. A Review 
of Transparent and Artificial Boundary Conditions Techniques for 
Linear and Nonlinear Schr\"odinger Equations. Comm. Comput. Phys. 4, 
729-796 (2008).

%

\bibitem{1}
X. Antoine and C. Besse. Unconditionally stable discretization 
schemes of non-reflecting boundary conditions for the one-dimensional 
Schrödinger equation. J. Comput. Phys., 188(1), 157-175 (2003).

%


\bibitem{bailung} H. Bailung, S. K. Sharma, and Y. Nakamura, Observation of Peregrine solitons in a multicomponent plasma with negative ions, Phys. Rev. Lett. 107, 255005 (2011).

\bibitem{6}
V.A. Baskakov and A.V. Popov. Implementation of transparent 
boundaries for numerical solution of the Schrödinger equation. Wave 
Motion, 14(2), 123-128 (1991).

%
\bibitem{Berenger}
J. B\'erenger. A perfectly matched layer for the absorption of electromagnetic waves. J. Comput. Phys. 114, 185-200, 1994.

\bibitem{anne}A. Boutet de Monvel, A.S. Fokas and D. Shepelsky, 
Analysis of the global relation for the nonlinear Schr\"odinger 
equation on the half-line, Lett. Math. Phys. 65, 199-212 (2003).

\bibitem{28}L.~Burgnies, O.~Vanb\'esien and D.~Lippens, Transient analysis of ballistic transport in stublike
quantum waveguides, Appl. Phys. Lett. 71, 803-805 (1997).

\bibitem{cha1}A. Chabchoub, N. P. Hoffmann, and N. Akhmediev, Rogue wave observation in a water wave tank, Phys. Rev. Lett. 106, 204502 (2011).

\bibitem{cha2} A. Chabchoub, N. Hoffmann, M. Onorato, and N. Akhmediev, Super rogue waves: observation of a higher-order breather in water waves, Phys. Rev. X 2, 011015 (2012).

\bibitem{cha3}A. Chabchoub, N. Hoffmann, H. Branger, C. Kharif, and 
N. Akhmediev, Experiments on wind-perturbed rogue wave hydrodynamics 
using the Peregrine breather model, Physics of Fluids 25 (2013) DOI: 10.1063/1.4824706

\bibitem{schober}A. Calini and C. M. Schober, Nat. Hazards Earth 
Syst. Sci., 14, 1431-1440 (2014).

\bibitem{33}J.F.~Claerbout, Coarse grid calculation of waves in inhomogeneous media with application to delineation of complicated seismic structure, Geophysics 35, 407-418 (1970).

%

\bibitem{dubard} P.~Dubard,  P.~Gaillard,  C.~Klein and V.B.~Matveev,   On multi-rogue wave solutions of the NLS 
equation and positon solutions of the KdV equation, Eur. Phys. J. 
Special Topics  185, 247-258 (2010)	


\bibitem{duque}J.~Duque, Solving time-dependent equations of 
Schr\"odinger-type using mapped infinite elements, Int. J. Mod. Phys. 
C \textbf{16}, No. 2 309-316 (2005).
 

\bibitem{hagstrom}T.~Hagstrom, Radiation boundary conditions for the 
numerical simulation of waves, Acta Numerica 8, 47-106 (1999). 

\bibitem{note}
S. G. Johnson. Notes on Perfectly Matched Layers (PMLs) (2010), http://math.mit.edu/~stevenj/18.369/pml.pdf.

%

\bibitem{kibler}B. Kibler, J. Fatome, C. Finot, G. Millot, F. Dias, 
G. Genty, N. Akhmediev, and J. M. Dudley, The Peregrine soliton in 
nonlinear fibre optics, Nat. Phys. 6, 790-795 (2010). 

\bibitem{etna} C.~Klein,  Fourth order time-stepping for low dispersion Korteweg-de 
Vries and nonlinear Schr\"odinger equation,  ETNA Vol. 29 116-135 (2008).

\bibitem{klein}
C. Klein and K. Roidot,  Fourth order time-stepping for 
Kadomtsev-Petviashvili and Davey-Stewartson equations, SIAM J. Sci. 
Comput., 33(6), 3333-3356. DOI: 10.1137/100816663 (2011). 


     \bibitem{KP2013} C.~Klein and R.~Peter, \emph{Numerical study of blow-up in solutions 
to generalized Korteweg-de Vries equations}, arXiv:1307.0603 	

\bibitem{pauline}X.~Antoine, C.~Besse, and P.~Klein. 
Absorbing Boundary Conditions for the One-Dimensional Schr\"odinger Equation with an Exterior Repulsive Potential, 
Journal of Computational Physics, 228(2), 312-335 (2009). 


\bibitem{ladouceur}F.~Ladouceur, Boundaryless beam propagation, Opt. 
Lett. 21, 4-5 (1996).

\bibitem{tau} C. Lanczos, Trigonometric interpolation of empirical 
and analytic functions, J. Math. and Physics, 17, 123-199 (1938)


\bibitem{87} M.F.~Levy, Parabolic equation models for electromagnetic wave propagation, IEE Electromagnetic
Waves Series 45 (2000).

\bibitem{95}C.W.~McCurdy, D.A.~Homer and T.N.~Resigno, Time dependent 
approach to collisional ionization
using exterior complex scaling, Phys. Rev. A 65, 042714 (2002).

%
%
\bibitem{nissan}
A. Nissen, G. Kreiss. An Optimized Perfectly Matched Layer for the 
Schrödinger Equation. Rapport technique, Department of Information 
Technology, Uppsala University (2009).

\bibitem{orszag}C.E. Grosch and S.A. Orszag, Numerical solution of 
problems in unbounded regions: coordinate transforms, J. Comput. 
Phys. 25, 273-296 (1977). 

\bibitem{Peregrine}D.H.~Peregrine, Water waves, nonlinear 
Schr\"odinger equations and their solutions, J. Austral. Math. Soc. 
B 25 16-43 (1983), doi:10.1017/S0334270000003891

\bibitem{112} F.~Schmidt and P.~Deuflhard, Discrete transparent 
boundary conditions for the numerical solution
of FresnelÕs equation, Comput. Math. Appl. 29, 53-76 (1995).

\bibitem{129} F.D. Tappert, The parabolic approximation method, in Wave Propagation and Underwater Acoustics, Lecture Notes in Physics 70, eds. J.B. Keller and J.S. Papadakis, Springer, New York,
  224-287 (1977).

\bibitem{semyon}
S.V.~Tsynkov. Numerical solution of problems on unbounded domains. A 
review. Appl. Numer. Math., 27(4), 465-532 (1998). 


\bibitem{trefethen} L.N. Trefethen, Spectral Methods in
    Matlab. SIAM, Philadelphia, PA (2000)


\bibitem{trefethenweb} www.comlab.ox.ac.uk/oucl/work/nick.trefethen

\bibitem{Lorene}www.lorene.obspm.fr

\bibitem{ZS} V.E. Zakharov and A.B.  Shabat, Exact theory of two-dimensional self-focusing and one-dimensional self-
modulation of waves in nonlinear media. Sov. Phys. JETP 34(1), 62-69 
(1972); translated from Zh. Eksp. Teor. Fiz. 1, 118-134 (1971) 

\bibitem{Zhengpml}
C. Zheng. A perfectly matched layer approach to the nonlinear 
Schr\"odinger wave equations. J. Comput. Phys. 227, 537-556 (2007).

\bibitem{Zhengtbc}
C. Zheng, Exact nonreflecting boundary conditions for one-dimensional 
cubic nonlinear Schr\"odinger equations, J. Comput. Phys. 
215, 552-565 (2006).


%
%
%
%
%
%

\end{thebibliography}
\end{document}